\documentclass[aos,MSNbibl,seceqn,dvips]{arximspdf}

\doi{10.1214/12-AOS1022} 
\volume{40}
\issue{3}
\pubyear{2012}
\firstpage{1737}
\lastpage{1763}

\makeatletter

\setattribute{abstract}{width}{280pt}

\newcommand{\rrvert}{\vert}
\newcommand{\llvert}{\vert}

\newcommand{\eps}{\varepsilon}
\newcommand{\rd}{{ d}}
\newcommand{\bR}{{\mathbb R}}
\newcommand{\bx}{{\mathbf{x}}}
\newcommand{\bv}{{\mathbf{v}}}
\newcommand{\bw}{{\mathbf{w}}}
\newcommand{\al}{\alpha}

\newcommand{\cal}{\mathcal}
\newcommand{\mG}{\mathcal G}

\newcommand{\st}{\dvtx}
\renewcommand{\epsilon}{\varepsilon}

\renewcommand{\P}{\mathbb{P}}
\newcommand{\E}{\mathbb{E}}
\newcommand{\R}{\mathbb{R}}
\newcommand{\C}{\mathbb{C}}
\newcommand{\N}{\mathbb{N}}

\newcommand{\tr}{\operatorname{Tr}}

\newcommand{\im}{\Im}

\newtheorem{theorem}{Theorem}[section]
\newtheorem{lemma}[theorem]{Lemma}
\newtheorem{corollary}[theorem]{Corollary}

\newproclaim{definition}[theorem]{Definition}
\newproclaim{example}[theorem]{Example}
\newproclaim{remark}[theorem]{Remark}

\newcommand{\mw}{m_W}
\newcommand{\I}{\mathcal{I}}
\newcommand{\bk}{\mathbf{k}}
\renewcommand{\b}{m}

\makeatother

\begin{document}
\begin{frontmatter}

\title{Edge universality of correlation matrices}
\runtitle{Edge universality of correlation matrices}

\begin{aug}
\author[A]{\fnms{Natesh S.} \snm{Pillai}\thanksref{t1}\ead[label=e1]{pillai@stat.harvard.edu}}
\and
\author[B]{\fnms{Jun} \snm{Yin}\corref{}\thanksref{t2}\ead[label=e2]{jyin@math.wisc.edu}}
\runauthor{N. S. Pillai and J. Yin}
\affiliation{Harvard University and University of Wisconsin--Madison}
\address[A]{Department of Statistics\\
Harvard University\\
1 Oxford Street\\
Cambridge, Massachusetts 02138\\
USA\\
\printead{e1}}
\address[B]{Department of Mathematics\\
University of Wisconsin--Madison\\
480 Lincoln Dr.\\
Madison, Wisconsin 53706\\
USA\\
\printead{e2}} 
\end{aug}

\thankstext{t1}{Supported by NSF Grant DMS-11-07070.}

\thankstext{t2}{Supported by NSF Grant DMS-10-01655.}

\received{\smonth{2} \syear{2012}}
\revised{\smonth{5} \syear{2012}}

%
\begin{abstract}
Let $\widetilde X_{M\times N}$ be a rectangular data matrix with independent
real-valued entries $[\widetilde x_{ij}]$ satisfying $\E\widetilde
{x}_{ij} = 0$
and $\E\widetilde{x}^2_{ij} = \frac{1 }{ M}$, $N, M\to\infty$. These
entries have a subexponential decay at the tails. We will be working
in the regime $N/M = d_N, \lim_{N\to\infty}d_N \neq0,1,\infty$. In
this paper we prove the edge universality of correlation matrices
${X}^\dagger X$, where the rectangular matrix $X$ (called the standardized
matrix) is obtained by normalizing each column of the data matrix
$\widetilde
X$ by its Euclidean norm. Our main result states that asymptotically
the $k$-point ($k \geq1$) correlation functions of the extreme
eigenvalues (at both edges of the spectrum) of the correlation matrix
${X}^\dagger X$ converge to those of the Gaussian correlation matrix, that
is, Tracy--Widom law, and,
thus, in particular, the largest and the smallest eigenvalues of
${X}^\dagger X$ after appropriate centering and rescaling converge to the
Tracy--Widom distribution. The asymptotic distribution of extreme
eigenvalues of the Gaussian correlation matrix has been worked out only
recently. As a corollary of the main result in this paper, we also
obtain that the extreme eigenvalues of Gaussian correlation matrices
are asymptotically distributed according to the Tracy--Widom law. The
proof is based on the comparison of Green functions, but the key
obstacle to be surmounted is the strong dependence of the entries of
the correlation matrix. We achieve this via a~novel argument which
involves comparing the moments of product of the entries of the
standardized data matrix to those of the raw data matrix. Our proof
strategy may be extended for proving the edge universality of other
random matrix ensembles with dependent entries and hence is of
independent interest.
\end{abstract}

%

\end{frontmatter}

\section{Introduction}\label{sec1} \label{secintro}
The aim\vspace*{1pt} of this paper is to prove the edge universality of correlation matrices.
The data matrix $\widetilde X = (\widetilde x_{ij}) $ is an $M \times
N$ matrix 
with independent centered\vadjust{\goodbreak} real-valued entries. The entries in each
column $j$ all are assumed to be identically distributed:
%
%
\begin{equation}
\label{eqnXmat} \widetilde x_{ij} = M^{-1/2}
q_{ij},\qquad \E q_{ij} = 0,\qquad \E q_{ij}^2 =
\sigma^2_j,\qquad 1 \leq i \leq M.
\end{equation}
Furthermore, the entries $q_{ij}$ have a subexponential decay,
that is, there exists
a constant $\vartheta>0$ such that for $u>1$,
%
%
\begin{equation}
\label{eqnXmatexpbd} \P\bigl(|q_{ij}| > u\sigma_j\bigr) \leq
\vartheta^{-1} \exp\bigl(-u^\vartheta\bigr).
\end{equation}
We will be working the regime
%
%
\begin{equation}
\label{eqnd} d=d_N=N/M,\qquad \lim_{ N\to\infty} d \neq0,1,\infty.
\end{equation}
Thus, without loss of generality, henceforth we will assume that for
some small constant $\theta$, for all $N \in\mathbb{N}$,
\[
\theta< d_N<\theta^{-1} \quad\mbox{and}\quad\theta<
|d_N-1|.
\]
Notice that all our constants may depend on $\theta$ and $\vartheta$,
but we will subsume
this dependence in the notation.

For a Euclidean vector $a \in\mathbb{R}^M$, define
the $\ell_2$ norm
\[
\| a \|_2:= \Biggl(\sum_{i=1}^M
a_i^2 \Biggr)^{1/2}.
\]
The matrix ${\widetilde{X}}^\dagger\widetilde X$ is the usual
covariance matrix.
The $j$th column of $\widetilde X$ is denoted by $\widetilde
\bx_j$.
Define the matrix $M \times N$ matrix $ X = ( x_{ij}) $
%
%
\begin{equation}
\label{defxnc} x_{ij}
:=\widetilde x_{ij}/\|\widetilde
\bx_j\|_2. 
\end{equation}
The $(N \times N)$ matrix ${X}^\dagger X$ is called the correlation
matrix.\setcounter{footnote}{2}\footnote{Some authors prefer to call this the standardized covariance
matrix, but we chose this terminology
from the statistical literature~\cite{John01}.}
Using the identity $\E x^2_{ij}=\frac1M \E\sum_ix^2_{ij}$, we have
\[
\E x^2_{ij}=M^{-1}.
\]
Since we are mainly interested in correlation matrices, without loss of
generality, henceforth we will assume that
\[
\sigma^2_j = 1 ,\qquad 1\leq j\leq N.
\]
Covariance matrices are ubiquitous in modern multivariate statistics
where the advance of technology has led to a profusion of
high-dimensional data sets. See~\cite{John01,John07,John08,NYcov} and
the references therein for motivation and applications in a
wide variety of fields.
Correlation matrices are sometimes preferred in certain statistical
applications. For instance, the classic exploratory method Principal
Component Analysis (PCA) is not invariant to change of scale in the
matrix entries. Therefore, it is often recommended first to standardize
the matrix entries and then perform PCA on the resulting correlation
matrix~\cite{John01}.

Recent progress in random matrix theory has led to a wealth of
techniques for proving universality of various matrix ensembles (see
\cite
{EKYY11,EKYY12,ESY1,ESY2,ESY3,ESY4,EPRSY,ESYY,EYYBulkuni,EYYgenwig,EYYrigid,Joha12,KY1,KY2,TaoVu09,TaoVu10}
and the references therein). Here the word universality refers to the
phenomenon that the asymptotic distributions of various functionals of
covariance/correlation matrices (such as eigenvalues, eigenvector,
etc.) are identical to those Gaussian covariance/correlation matrices.
Thus, harnessing these methods to obtain universality results in
statistical problems is an important step, since these results let us
calculate the exact asymptotic distributions of various test statistics
without having restrictive distributional assumptions of the matrix
entries. For instance, an important consequence of universality is that
in some cases one can perform various hypothesis tests under the
assumption that the matrix entries are \textit{not} normally distributed
but use the same test statistic as in the Gaussian case.

In this context, in a recent paper~\cite{NYcov} we studied the
asymptotic distribution of the eigenvalues of the covariance matrix
${\widetilde{X}}^\dagger\widetilde X$ under the assumptions of
(\ref{eqnXmat}) and~(\ref{eqnXmatexpbd}). In~\cite{NYcov}, we proved
that the Stieltjes transform of the empirical eigenvalue distribution
of the sample covariance matrix is given by the Marcenko--Pastur law
\cite{MP} uniformly up to the edges of the spectrum with an error of
order $ (N \eta)^{-1}$, where $\eta$ is the imaginary part of the
spectral parameter in the Stieltjes transform. From this strong local
Marcenko--Pastur law, we derived the following results: (1)~rigidity of
eigenvalues (2) delocalization of eigenvectors (3) universality of
eigenvalues in the bulk and (4) universality of eigenvalues at the
edges. Furthermore, in our proof of edge universality of eigenvalues
for covariance matrices (see Theorem~7.5 of~\cite{NYcov}), we gave a
sufficient criterion for checking whether two matrices of form
$Q^\dagger Q$ ($Q$ is a data matrix) have the same asymptotic
eigenvalue distribution at the edge (see Section~\ref{secproofsketch}
for details). Here ${Q}^\dagger Q$ could be quite general, including
covariance and correlation matrices.

Verifying the above criteria for correlation matrices is much more
complicated, owing to the fact that even if it has the same form
${X}^\dagger X$ as above, the matrix entries of $X $ are not independent.
Fortunately in~\cite{NYcov}, as a byproduct, we also proved the strong
Marcenko--Pastur law, the rigidity of eigenvalues and delocalization of
eigenvectors of correlation matrices (see Lemma~\ref{lemstrMplawri}
in Section~\ref{sec2} below or Theorem 1.5 of~\cite{NYcov}). In this
paper, we
complete the research program initiated in~\cite{NYcov} by proving the
edge universality of correlation matrices. There are not many papers
which study the asymptotics of the correlation matrices as compared to
the relatively large literature on covariance matrices. The asymptotic
distribution of the largest (appropriately rescaled) eigenvalue of the
Gaussian correlation matrix was only very recently established by
\cite{Bao1}. As will be explained below, we also obtain this result as a
special case of our main result and, more importantly, \textit{we do not
need this result in our proof} (see Remark~\ref{remGLaw}). The almost sure
convergence of the largest and smallest eigenvalues of the correlation
matrix was established in~\cite{Jiang}. The very recent paper
\cite{Bao1}, relying on our results in~\cite{NYcov}, shows that the
asymptotic distribution of the largest or smallest eigenvalue of the
correlation matrix is given by the Tracy--Widom law, under the
assumption that the data matrix $X$ satisfies~(\ref{eqnXmat}) and its
entries have \textit{symmetric} distributions. In particular, the
authors in~\cite{Bao1} use the above mentioned sufficiency criteria
for edge universality developed in~\cite{NYcov}. Furthermore, the
assumption that the matrix entries are symmetric is very restrictive
and not natural in statistical applications. In this paper we will
build on our previous work~\cite{NYcov} and prove edge universality of
correlation matrices just under the assumptions~(\ref{eqnXmat}) and
(\ref{eqnXmatexpbd}). Furthermore, we believe that all of our main
results should hold if one replaces
the subexponential tail decay of the matrix entries by a uniform bound
on the $p$th moment $(p > 4)$ of the matrix entries (e.g.,
$p=13$ will suffice), as proved in~\cite{EKYY12} for Wigner matrices.

The central ideas in this paper are based on the general machinery for
proving universality established in a series of recent papers
\cite
{EKYY11,EKYY12,ESY1,ESY2,ESY3,ESY4,EPRSY,ESYY,EYYBulkuni,EYYgenwig,EYYrigid,KY1,KY2},
where the authors
Yau, Erd{\H o}s et al. study the distribution of eigenvalues
and eigenvectors by studying the Green's functions (resolvent) of the
random matrices.

The proof of this paper is based on the comparison of Green's functions
first initiated in~\cite{EYYBulkuni}, but, as mentioned earlier, the
key obstacle to be surmounted is the strong dependence of the entries
of the correlation matrix. We achieve this via a novel argument which
involves comparing the moments of the product of the entries of the
standardized data matrix to those of the raw data matrix (see
Section~\ref{secproofsketch} for a summary of the key ideas). Our proof
strategy may be extended for proving the edge universality of other
random matrix ensembles with dependent entries and hence is of
independent interest. Furthermore, it will be interesting to see if
bulk universality of correlation matrices can be established using the
methods developed in this paper.

Let us state the main result now. We denote $\lambda_i$, $1\leq i\leq
N$, as the eigenvalues of $X^\dagger X$ and $\lambda_\alpha=0$ for $
\min\{N, M\}+1\leq\alpha\leq\max\{N,M\} $. We order them as
\[
\lambda_1\geq\lambda_2\geq\cdots\geq\lambda_{\max\{M,N\}}
\geq0.
\]
Analogously, let $\widetilde{\lambda}_\alpha$ denote the eigenvalues
values of the matrix ${\widetilde X}^\dagger{\widetilde X}$.

The following is the main result of this paper. It shows that the
largest and smallest $k$ eigenvalues of the correlation matrix, after
appropriate centering and rescaling, converge in distribution to those
of the corresponding covariance matrix.
%
%
\begin{theorem}[(Edge universality)] \label{thmmain}
Let $X$ and $\widetilde X$, respectively, denote the correlation and covariance
matrix as defined in~(\ref{eqnXmat})--(\ref{defxnc}). For any fixed
$k\in\mathbb{N}$,
there exists ${\varepsilon}> 0$ and $\delta>0$ such that for any $\{
s_1, s_2, \ldots, s_k\} \in\mathbb{R}$ (which may depend on~$N$),
there exists $N_0 \in\mathbb{N}$ independent of $s_1, s_2, \ldots,
s_k$ such that for all $N \geq N_0$, we have
%
%
\begin{eqnarray}
&&
\P\bigl( N^{2/3} ( \widetilde\lambda_1 - \lambda_+)
\leq s_1- N^{-{\varepsilon}}, \ldots, N^{2/3} ( \widetilde
\lambda_k - \lambda_+) \leq s_{k }- N^{-{\varepsilon}}
\bigr)- N^{-\delta}
\nonumber\hspace*{-30pt}
\\
&&\qquad \leq\P\bigl( N^{2/3} ( \lambda_1 - \lambda_+) \leq
s_1, \ldots, N^{2/3} ( \lambda_k - \lambda_+)
\leq s_{k } \bigr)\hspace*{-30pt}
\nonumber\\[-8pt]\\[-8pt]
&&\qquad \leq\P\bigl( N^{2/3} ( \widetilde\lambda_1 -
\lambda_+) \leq s_1+ N^{-{\varepsilon}}, \ldots, N^{2/3} (
\widetilde\lambda_k - \lambda_+) \leq s_{k }+
N^{-{\varepsilon}} \bigr)\nonumber\hspace*{-30pt}\\
&&\qquad\quad{}+ N^{-\delta}.
\nonumber\hspace*{-30pt}
\end{eqnarray}
An analogous result holds for the $k$ smallest eigenvalues. 
\end{theorem}

In~\cite{Pech07,Sod1} and~\cite{Sosh1}, Peche, Soshnikov
and Sodin proved that for some covariance matrices (including the
Wishart matrix), the largest and smallest $k$ eigenvalues after
appropriate centering\vspace*{1pt} and rescaling converge in distribution to the
Tracy--Widom law\footnote{Here we use the term Tracy--Widom law as in
\cite{Sosh1}.} whose density is a smooth function. Combining with our
recent result on the universality of covariance matrices in
\cite{NYcov}, we have the following immediate corollary for Theorem \ref
{thmmain}:
%
%
\begin{corollary} Let $X$ denote the correlation matrix as defined
in~(\ref{eqnXmat})--(\ref{defxnc}). For any fixed $k>0$, we have
\begin{eqnarray*}
&&
\biggl(\frac{M\lambda_1 - (\sqrt{N} + \sqrt{M})^2} {(\sqrt{N} +
\sqrt{M})({1 }/{\sqrt N} +
{1}/{\sqrt M})^{1/3}}, \ldots,\\
&&\hspace*{8pt}\qquad\frac{M\lambda_k- (\sqrt{N} +
\sqrt{M})^2} {(\sqrt{N} + \sqrt{M})({1 }/{\sqrt N} +
{1}/{\sqrt M})^{1/3}} \biggr) \\
&&\qquad\longrightarrow
\mathrm{TW}_1,
\end{eqnarray*}
where $\mathrm{TW}_1$ denotes the Tracy--Widom distribution. An analogous
statement holds for the $k$-smallest (nontrivial) eigenvalues.
\end{corollary}
%
%
\begin{remark} \label{remGLaw}
Thus, as a special case, we also obtain the TW law for the Gaussian
correlation matrices.
\end{remark}
Although the current paper builds on our recent work~\cite{NYcov}, it
is mostly self-contained and for the reader's convenience, we will
recall all of the needed results from~\cite{NYcov}. The rest of the
paper is organized as follows. In Section~\ref{secprelim}, after
establishing some notation, we give the key results establishing the
strong Marcenko--Pastur law and rigidity of eigenvalues for correlation
matrices,\vadjust{\goodbreak} as obtained from~\cite{NYcov}.
In Section~\ref{secproofsketch} we give a brief proof sketch
illustrating the key ideas.
In Section~\ref{secpromainre} we give the proof of the main results
and in Section~\ref{secmomcomp} we prove some technical lemmas which
constitute the key ingredients in the proof of the main result. For the
rest of the paper the letter $C$ will denote a generic constant whose
value might change from one line to the next, but will be independent
of everything else.
The notation $O_\epsilon(N^a)$ will be used to denote $O(N^{a +
C{\varepsilon}})$.

\section{Preliminaries}\label{sec2} \label{secprelim} We will adopt the notation
used in this paper from~\cite{NYcov}.
Define the Green function of ${X}^\dagger X$ by
%
%
\begin{equation}
\label{eqngreen} G_{ij}(z) = \biggl(\frac1{{X}^\dagger X-z}
\biggr)_{ij},\qquad z=E+i\eta,\qquad E\in\bR, \eta>0.
\end{equation}
The Stieltjes transform of the empirical
eigenvalue distribution of ${X}^\dagger X $ is given by
%
%
\begin{equation}\label{mNdef}
m(z) := \frac{1}{N} \sum_j
G_{jj}(z) = \frac{1}{N} \tr\frac
{1}{{X}^\dagger X-z}.
\end{equation}
%
Recall that $d = N/M$ from~(\ref{eqnd}) and
define
%
%
\begin{equation}
\label{eqnlambdapm} \lambda_\pm:= \bigl(1\pm\sqrt{d} \bigr)^2.
\end{equation}
The Marcenko--Pastur (henceforth abbreviated by MP) law is given by
%
%
\begin{equation}
\label{defrhow} \varrho_W(x)=\frac{1}{2\pi d}\sqrt{
\frac{ [(\lambda_+-x)(x-\lambda_-) ]_+}{x^2}}.
\end{equation}
We define $m_W (z)$, $z\in\C$, as the Stieltjes transform of $\varrho
_W$, that is,
%
%
\begin{equation}
\label{defmW} m_W (z)=\int_{\R}
\frac{\varrho_W(x)}{(x-z)} \,dx.
\end{equation}
The function
$m_W$ depends on $d$ and has the closed form solution
%
%
\begin{equation}
\label{mwz=} m_W(z)=\frac{1-d-z +i \sqrt{(z-\lambda_-)(\lambda_+-z)}}{2
\,d z},
\end{equation}
where $\sqrt{\ }$ denotes the square root on a complex
plane whose branch cut is the negative real line.
We also define the classical location of the eigenvalues with $\rho_W$
as follows:
%
%
\begin{equation}
\label{defgj} \int_{\gamma_j}^{\lambda_+}\varrho_W(x)
\,dx =\int_{\gamma
_j}^{+\infty}\varrho_W(x)
\,dx=j/N.
\end{equation}
%
%
Define the parameter
%
%
\begin{equation}
\label{eqnvarphi} \varphi:= (\log N)^{\log\log N}.
\end{equation}

%
%
\begin{definition}[(High probability events)]\label{defhp}
Let $\zeta> 0$.
We say that an 
event $\Omega$ holds with \textit{$\zeta$-high probability} if there\vadjust{\goodbreak}
exists a constant $C > 0$ such that
%
%
\begin{equation}
\label{highprob} \P\bigl(\Omega^c\bigr) \leq N^C \exp
\bigl(-\varphi^\zeta\bigr)
\end{equation}
for large enough $N$.
\end{definition}

Let us first give the following large deviation lemma for independent
random variables (see~\cite{EYYBulkuni}, Appendix B for a proof).
%
%
\begin{lemma}[(Large deviation lemma)] \label{lemlargdev}
Suppose, for $1 \leq i \leq M$, $a_i$ are independent, mean $0$ complex
variables, with
$\E|a_i|^2 = \sigma^2$ and have a subexponential decay as in (\ref
{eqnXmatexpbd}). Then there exists a constant $\rho\equiv\rho
(\vartheta) > 1$ such that, for any $\zeta> 0$ and for any $A_i \in
\mathbb{C}$ and $B_{ij} \in\mathbb{C}$, the bounds
%
%
\begin{eqnarray}
\label{eqnld1}
\sum_{i=1}^M a_i
A_i &\leq& (\log M)^{\rho\zeta} \sigma\|A\|,
\\
\label{eqnld2}
\Biggl|\sum_{i=1}^M \bar a_i
B_{ii} a_i - \sum_{i=1}^M
\sigma^2 B_{ii} \Biggr| &\leq& (\log M)^{\rho\zeta}
\sigma^2 \Biggl(\sum_{i=1}^M
|B_{ii}|^2\Biggr)^{1/2},
\\
\label{eqnld3}
\biggl|\sum_{i \neq j} \bar a_i
B_{ij} a_j\biggr| &\leq& (\log M)^{\rho\zeta}
\sigma^2 \biggl(\sum_{i \neq j}
|B_{ij}|^2\biggr)^{1/2}
\end{eqnarray}
hold with $\zeta$-high probability.
\end{lemma}

It can be easily seen that for any fixed $j \leq N$, the random variables
defined by $a_i = x_{ij}, 1 \leq i \leq M$, satisfy the large deviation
bounds~(\ref{eqnld1}),~(\ref{eqnld2})
and~(\ref{eqnld3}), for any $A_i \in\mathbb{C}$ and $B_{ij} \in
\mathbb{C}$ and $\zeta> 0$.

Thus, the main result of~\cite{NYcov} (see Theorem 1.5 of
\cite{NYcov}) is applicable for the correlation matrix $X$, yielding the
following strong local MP law and rigidity of eigenvalues:
%
%
\begin{lemma}[(Strong local Marcenko--Pastur law and rigidity of the
eigenvalues of the correlation matrix)] \label{lemstrMplawri}
Let $X = [x_{ij}]$ be the correlation matrix given by~(\ref{defxnc}).
Then for any $\zeta>0$ there exists a constant $C_\zeta$ such that
the following events hold with $\zeta$-high probability.

\begin{longlist}
\item
The Stieltjes transform of the empirical
eigenvalue distribution of $X^\dagger X$ satisfies
%
%
\begin{equation}
\label{Lambdafinal} \bigcap_{z \in{ \b S}(C_\zeta)} \biggl\{{\bigl
\llvert m(z)-m_W(z)\bigr\rrvert\leq\varphi^{C_\zeta}
\frac{1}{N\eta}} \biggr\},
\end{equation}
where $ { \b S}(C_\zeta) $ defined as the set
\[
{ \b S}(C_\zeta):= \bigl\{z \in\C\st{\mathbf1}_{d>1}(
\lambda_-/5)\leq E \leq5\lambda_+, \varphi^{C_\zeta} N^{-1} \leq
\eta\leq10(1+d) \bigr\}.\vadjust{\goodbreak}
\]

\item The individual matrix elements of
the Green function satisfy
%
%
\begin{equation}
\label{Lambdaofinal} \bigcap_{z \in{ \b S}(C_\zeta)} \Biggl\{{\bigl
\llvert G_{ij}(z)-m_W(z)\delta_{ij}\bigr
\rrvert\leq\varphi^{C_\zeta} \Biggl(\sqrt{\frac{\Im \mw(z) }{N\eta}}+
\frac{1}{N\eta} \Biggr)} \Biggr\}.
\end{equation}
\item The smallest nonzero and largest eigenvalues of $X^\dagger X$ satisfy
%
%
\begin{equation}
\label{443} \lambda_- -N^{-2/3}\varphi^{C_\zeta}\leq
\min_{j\leq\min\{M, N\}} \lambda_j \leq\max_{j }
\lambda_j \leq\lambda_++N^{-2/3}\varphi^{C_\zeta}.
\end{equation}
\item Rigidity\vspace*{1pt} of the eigenvalues: recall $\gamma_j$ in~(\ref{defgj}).
For any $1\leq j\leq\break\min\{M,N\}$,
let $ \widetilde j=\min\{\min\{N, M\}+1-j, j \}. $
Then
%
%
\begin{equation}
\label{resrig} \llvert\lambda_j-\gamma_j\rrvert
\leq\varphi^{C_\zeta} N^{-2/3} \widetilde j^{-1/3}.
\end{equation}
\end{longlist}
\end{lemma}

We conclude this section with the following theorem quoted from
\cite{NYcov} (see Theorem 1.7 in~\cite{NYcov}) on edge universality of
covariance matrices, which is also needed for our proof of the edge
universality of the correlation matrix. Define two independent matrices
$\widetilde{X}{}^{\bv} = [\widetilde {x}^{\bv}_{ij}], \widetilde
X^{\bw}= [\widetilde{x}^{\bw}_{ij}]$ with the entries
$\widetilde{x}^{\bv}_{ij},\widetilde{x}^{\bw}_{ij}$ satisfying (\ref
{eqnXmat}) and~(\ref{eqnXmatexpbd}) and the entries $\widetilde{x}^{\bv
}_{ij},\widetilde{x}^{\bw}_{ij}$ are mutually independent.
Henceforth,\vspace*{1pt} we will write $\mathbb{E}^{\bv}, \P^{\bv}$
($\mathbb{E}^{\bw}, \P^{\bw }$) to indicate that the expectation and
probability are computed for the ensemble $\widetilde{X}{}^{\bv}$,
($\widetilde{X}{}^{\bw}$).
%
%
\begin{theorem}[(Universality of extreme eigenvalues of covariance
matrices)] \label{twthm}
There exists ${\varepsilon}>0$ and $\delta>0$
such that
for any $s \in\mathbb{R}$ (which may depend on $N$) there exists $N_0
\in\mathbb{N}$ independent of $s$ such that for all $N \geq N_0$,
we have
%
%
\begin{eqnarray}
\label{tw}
&&
\P^\bv\bigl( N^{2/3} \bigl( \widetilde
\lambda_1^\bv- \lambda_+\bigr) \leq s- N^{-{\varepsilon}}
\bigr)- N^{-\delta} \nonumber\\
&&\qquad\leq \P^\bw\bigl( N^{2/3} \bigl(
\widetilde\lambda_1^\bw- \lambda_+\bigr) \leq s \bigr)
\\
&&\qquad \leq \P^\bv\bigl( N^{2/3} \bigl( \widetilde
\lambda_1^\bv- \lambda_+\bigr) \leq s+ N^{-{\varepsilon}}
\bigr)+ N^{-\delta}.
\nonumber
\end{eqnarray}
An analogous result holds for the smallest eigenvalues $\widetilde
\lambda^{\bv}_{\min\{M,N\}}$
and $\widetilde\lambda^{\bw}_{\min\{M,N\}}$.
\end{theorem}


As remarked in~\cite{NYcov}, Theorem~\ref{twthm} can be extended to
finite correlation functions of extreme eigenvalues as follows:
%
%
\begin{eqnarray}
\label{twa}
&& \P^\bv\bigl( N^{2/3} \bigl( \widetilde
\lambda_1^\bv- \lambda_+\bigr) \leq s_1-
N^{-{\varepsilon}}, \ldots, N^{2/3} \bigl( \widetilde\lambda_k^\bv-
\lambda_+\bigr) \leq s_{k }- N^{-{\varepsilon}} \bigr)- N^{-\delta}
\nonumber\hspace*{-30pt}
\\
&&\qquad \leq\P^\bw\bigl( N^{2/3} \bigl( \widetilde
\lambda_1^\bw- \lambda_+\bigr) \leq s_1,
\ldots, N^{2/3} \bigl( \widetilde\lambda_k^\bw-
\lambda_+\bigr) \leq s_{k } \bigr)\hspace*{-30pt}
\nonumber\\[-8pt]\\[-8pt]
&&\qquad \leq\P^\bv\bigl( N^{2/3} \bigl( \widetilde
\lambda_1^\bv- \lambda_+\bigr) \leq s_1+
N^{-{\varepsilon}}, \ldots, N^{2/3} \bigl( \widetilde\lambda_k^\bv-
\lambda_+\bigr) \leq s_{k }+ N^{-{\varepsilon}} \bigr)\nonumber\hspace*{-30pt}\\
&&\qquad\quad{}+ N^{-\delta}
\nonumber\hspace*{-30pt}
\end{eqnarray}
for all $k$ fixed and sufficiently large $N$. We remark that edge
universality is usually
formulated in\vadjust{\goodbreak} terms of joint distributions of edge eigenvalues
as in~(\ref{twa}) with fixed parameters $s_1, s_2,\ldots$
etc. However, we note that
Theorem~\ref{twthm}
holds uniformly in these parameters, and thus
they may depend on~$N$.

\section{Key ideas and proof sketch}\label{sec3} \label{secproofsketch}
Our basic strategy is the so-called ``Green function comparison'' method
initiated in a recent series of papers including
\mbox{\cite{EYYBulkuni,EYYgenwig,EYYrigid}} for proving universality
for (generalized) Wigner matrices. The Green function comparison method
has subsequently been applied to proving the spectral universality of
adjacency matrices of random graphs~\cite{EKYY11,EKYY12}, the
universality of eigenvectors of Wigner matrices~\cite{KY1}, as well as
the the spectrum of additive finite-rank deformations of Wigner
matrices and the isotropic local semicircle law~\cite{KY2}.\vspace*{1pt}


In this paper, we will show that~(\ref{tw}) and~(\ref{twa}) still
hold with $\widetilde X^\bv$ and~$\widetilde X^\bw$ replaced by the
correlation
matrix $X$
and the corresponding covariance matrix $\widetilde X$, that is,
Theorem \ref
{thmmain}. To show this result, we introduce a sufficient criteria for
(\ref{tw}) and~(\ref{twa}) derived in~\cite{NYcov} (see Theorem 7.5
of~\cite{NYcov}).

Consider two matrix ensembles $ X^{\bv}, X^{\bw}$ (could be
covariance, correlation or more general matrix\footnote{Notice that
throughout the paper we use $X$ for the correlation matrix and
$\widetilde X$
for the covariance matrix. This is the only instance we denote a
generic matrix by $X$ for compactness of notation.}) and let their
respective Green functions and empirical Stieltjes transforms [see
(\ref{eqngreen}) and~(\ref{mNdef})] be denoted by $G^{\bv}, G^{\bw
}$ and~$m^{\bv}, m^{\bw}$. To prove that the asymptotic distribution
of the
extreme eigenvalues of the matrix ensembles $ X^{\bv}, X^{\bw}$ are
identical in the sense of~(\ref{tw}) and~(\ref{twa}), it suffices to
show the following~\cite{NYcov}:
\begin{longlist}
\item The matrices $ X^{\bv}, X^{\bw}$ satisfy the strong
Marcenko--Pastur law and the rigidity of eigenvalues as given in Lemma
\ref{lemstrMplawri}.
\item The difference of the expectation of smooth functionals of the
corresponding Green functions ($G^{\bv}, G^{\bw}$ and $m^{\bv},
m^{\bw}$) evaluated at the spectral edge must vanish asymptotically.
More precisely, as pointed out in~\cite{NYcov}, it suffices to
establish Theorems~\ref{GFCT} and~\ref{GFCT2} below for the matrices $
X^{\bv}, X^{\bw}$.
\end{longlist}

%
\begin{theorem}[(Green function comparison theorem on the edge)] \label{GFCT}
Let $F\dvtx\R\to\R$ be a function whose derivatives satisfy
%
%
\begin{equation}
\label{gflowder} \max_{x}\bigl|F^{(\al)}(x)\bigr| \bigl(|x|+1
\bigr)^{-C_1} \leq C_1,\qquad \al=1, 2, 3, 4,
\end{equation}
for some constant $C_1>0$.
Then there exist ${\varepsilon}_0>0, N_0 \in\mathbb{N}$ and $\delta
> 0$ depending only on $C_1$ such that for any ${\varepsilon
}<{\varepsilon}_0$, $N \geq N_0$ and real numbers $E$, $E_1$ and $E_2$
satisfying
\[
|E-\lambda_+|\leq N^{-2/3+\eps},\qquad |E_1-\lambda_+|\leq
N^{-2/3+\eps},\qquad |E_2-\lambda_+|\leq N^{-2/3+\eps}
\]
and $\eta_0 = N^{-2/3-{\varepsilon}}$,
we have
%
%
\begin{equation}
\label{maincomp}
\bigl|\E^\bv F \bigl( N \eta_0 \Im
m^{\bv} (z) \bigr) - \E^\bw F \bigl( N \eta_0
\Im m^{\bw} (z) \bigr) \bigr| \leq C N^{-\delta+C {\varepsilon}},\quad z=E+i
\eta_0,\hspace*{-35pt}
\end{equation}
and
%
%
\begin{eqnarray}
\label{c1}\qquad
&&
\biggl\llvert\E^\bv F \biggl(N \int_{E_1}^{E_2}
\rd y \,\Im m^{\bv} (y +i\eta_0) \biggr)-\E^\bw F
\biggl(N \int_{E_1}^{E_2} \rd y\, \Im m^{\bw}
(y+i\eta_0 ) \biggr)\biggr\rrvert\nonumber\\[-8pt]\\[-8pt]
&&\qquad\leq C N^{-\delta+C {\varepsilon}}\nonumber
\end{eqnarray}
for some constant $C$.
\end{theorem}

%
\begin{theorem} \label{GFCT2}
Fix any $k\in\N_+$ and
let $F\dvtx\R^k \to\R$ be a smooth, bounded function with bounded derivatives.
Then there exist ${\varepsilon}_0>0, N_0 \in\mathbb{N}$ and $\delta
> 0$ such that for any ${\varepsilon}<{\varepsilon}_0$, $N \geq N_0$
and sequence of real numbers $E_k < \cdots< E_{1}< E_0 $
with $|E_j-\lambda_+|\leq N^{-2/3+{\varepsilon}}$, $j=0,1,\ldots, k$
and $\eta_0 = N^{-2/3 - {\varepsilon}}$, we have
%
%
\begin{eqnarray}\label{c11}
&&
\biggl\llvert\E^{\bv} F \biggl(N \int_{E_1}^{E_0}
\rd y \,\Im m^{\bv} (y +i\eta_0), \ldots, N \int
_{E_k}^{E_0} \rd y \,\Im m^{\bv} (y+i
\eta_0 ) \biggr)\nonumber\\
&&\hspace*{188.3pt}{} - \E^{\bw} F \bigl(m^\bv
\rightarrow m^\bw\bigr) \biggr\rrvert\\
&&\qquad\leq N^{-\delta},\nonumber
\end{eqnarray}
where the second term in the left-hand side above is obtained by
changing the arguments of $F$ in the first term from $m^\bv$ to $m^\bw
$ and keeping all the other parameters fixed.
\end{theorem}

%
\begin{remark}
Theorems~\ref{GFCT} and~\ref{GFCT2} yield the edge universality of
the $k$-point correlation functions at the edge for $k=1$ and $k \geq
1$, respectively.
\end{remark}
Thus, to complete the proof of Theorem~\ref{thmmain}, by the Green
function comparison method
it suffices to show (i) and (ii) above for
\[
X^\bv=X,\qquad X^\bw=\widetilde X,
\]
where $X^\dagger X$ denotes the correlation matrix and $ \widetilde
X^\dagger
\widetilde X$ is the corresponding covariance matrix.
Here condition (i) is guaranteed by Theorem~\ref{lemstrMplawri}.

Verifying condition (ii) entails the heart of this paper.
In previous works mentioned earlier, the authors use a Lindeberg
replacement strategy, as in \cite{Chat,TaoVu09}. These proofs
proceed via showing that the distribution of some smooth functional of
the Green function (e.g., $G_{ii}$, $m$ and $\langle\bx_1,
G \bx_1\rangle$) of the two matrix ensembles is identical asymptotically
provided that the first two (in some cases up to four) moments of all
matrix elements of these two ensembles are identical. For instance, if
one needs to show the edge universality\vadjust{\goodbreak} of two covariance matrices
$\widetilde
{X}{}^{\bv}$ and $\widetilde{X}{}^{\bw}$, the basic strategy is to express
%
%
\begin{equation}
\label{eqnfineq0} \E F\bigl(\widetilde G^{\bv}\bigr) -\E F\bigl(
\widetilde G^{\bw}\bigr) = \sum_{\gamma
=1}^{MN}
\E F(\widetilde G_{\gamma}) -\E F(\widetilde G_{\gamma- 1}),
\end{equation}
where $F$ is a smooth function and $\widetilde G_{\gamma}$ denotes the
Green function of the ensemble $\widetilde X_{\gamma}$ (with
$\widetilde X_{0} = \widetilde X^{\bv}$) which\vspace*{1pt} is obtained
from $\widetilde X_{\gamma-1}$ by replacing the distribution of the
$ij$th entry of $\widetilde {X}_{\gamma -1}[{ij}]$ with $\widetilde
X^{\bw}[{ij}]$ [here $ \gamma= i+ (j-1)M$] so that $\widetilde X_{MN} =
\widetilde X^{\bw}$. The next step is to obtain an estimate
%
%
\begin{equation}
\label{gg1} \E F(\widetilde G_{\gamma}) -\E F(\widetilde
G_{\gamma- 1}) =o\bigl(N^{-2}\bigr) 
\end{equation}
for each of the $N^2$ terms in the sum~(\ref{eqnfineq0}).
Usually\vspace*{1pt}~(\ref{gg1}) is obtained by resolvent expansions,
perturbation theory and the fact that $\widetilde X_\gamma$ and
$\widetilde X_{\gamma-1}$ differ by a single entry and the first few
moments of these two distributions are the same.

But clearly the above method does not work in our case, since the
entries within the same column are not independent and, therefore, one
cannot replace the distribution of a single entry of a column without
changing the distribution of all the other $M-1$ entries. To circumvent
this, in~\cite{NYcov} a~new telescoping argument consisting of $O(N)$
ensembles was used for the comparison of Green functions. The idea is
that instead of replacing entries one at a time, one can replace the
entries of the data matrix column by column and thus require only
$O(N)$ ensembles. This argument from~\cite{NYcov} is adapted here
along with new insights for dealing with nonindependence of the entries
and is outlined below.

Now we set $X^\bv=X, X^\bw=\widetilde X$.
For $1\leq\gamma\leq N$, let $X_\gamma$ denote the random matrix
whose $j$th column is the same as that of $X^{\bv}$ if $j > \gamma$
and that of $X^{\bw}$
otherwise. In particular, we can choose $X_0 = X^{ \bv} = X$ and
$X_{N} = X^{\bw} = \widetilde X$, where $X$ is correlation matrix and
$\widetilde X$
the corresponding covariance matrix of $X$. As before, we define
\[
m_\gamma(z) =\frac{1}{N}\tr G_\gamma(z),\qquad
G_\gamma(z)=\bigl(X_\gamma^\dagger X_\gamma-z
\bigr)^{-1},
\]
so that we have telescoping sum
%
%
\begin{eqnarray}\label{eqntelesum}
&&
\E^\bw F \bigl( N \eta_0 \Im m^{\bw} (z)
\bigr) - \E^\bv F \bigl( N \eta_0 \Im m^{\bv}
(z) \bigr) \nonumber\\[-8pt]\\[-8pt]
&&\qquad= \sum_{\gamma= 1}^N\E F \bigl( N
\eta_0 \Im m_{\gamma} (z) \bigr) - \E F \bigl( N
\eta_0 \Im m_{\gamma-1} (z) \bigr).\nonumber
\end{eqnarray}
Clearly,~(\ref{maincomp}) will follow from~(\ref{eqntelesum}) and
the following estimate:
%
%
\begin{equation}
\label{maincomp20} \bigl|\E F \bigl( N \eta_0 \Im m_{\gamma} (z)
\bigr) - \E F \bigl( N \eta_0 \Im m_{\gamma-1} (z) \bigr) \bigr| \leq
O_{\varepsilon}\bigl(N^{-1-\delta}\bigr)
\end{equation}
for some $\delta> 0$.
Our strategy to obtain~(\ref{maincomp20}) is the following.
First notice that
\begin{eqnarray*}
&&
\E F \bigl( N \eta_0 \Im m_{\gamma} (z) \bigr) - \E F \bigl(
N \eta_0 \Im m_{\gamma-1} (z) \bigr)\\
&&\qquad = \E F \bigl(
\eta_0 \im\tr G_\gamma(z) \bigr) - \E F \bigl(
\eta_0 \im\tr G_{\gamma-1} (z) \bigr).
\end{eqnarray*}
%
Let $X^{(\gamma)}$ be the $M \times(N-1)$ matrix obtained by removing
the $\gamma$th column of $X_\gamma$, which has the same
distribution of the $M \times(N-1)$ matrix obtained by removing the
$\gamma$th column of $X_{\gamma-1}$. Define
%
%
\begin{equation}
\label{dwbb} G^{(\gamma)}= \bigl(\bigl(X^{(\gamma)}\bigr)^\dagger
\bigl(X^{(\gamma)}\bigr)-z \bigr)^{-1},\qquad \mu= \eta_0
\im\tr G^{(\gamma)} - \im\frac{\eta_0}{z}.
\end{equation}
In Lemma~\ref{premainlemma} we will establish~(\ref{maincomp20}) by
showing that
%
%
\begin{eqnarray}
\label{ss440}
&&\bigl(\E F ( \eta_0 \im\tr G_\gamma) - \E F
( \mu) \bigr) - \bigl(\E F ( \eta_0 \im\tr G_{\gamma-1} ) - \E
F ( \mu) \bigr) \nonumber\\[-8pt]\\[-8pt]
&&\qquad= O_{\varepsilon}\bigl(N^{-7/6}\bigr).\nonumber
\end{eqnarray}
Once~(\ref{maincomp20}) is verified, the main result follows by virtue
of Theorems~\ref{GFCT} and~\ref{GFCT2} as mentioned in the beginning
of this section. Notice that since the columns of the data matrix $
X^\bv$, $X^\bw$ are assumed to be independent, $\mu$ is independent
of the $\gamma$th column of $ X^\bv$, $X^\bw$ or,
equivalently, the $\gamma$th column of $ X_\gamma$,
$X_{\gamma-1}$.

Thus, it boils down to establishing~(\ref{ss440}) in the case $X_0 =
X^{ \bv} = X$ and $X_{N} = X^{\bw} = \widetilde X$.
Our proof relies on the key observation that even if the entries of the
$\gamma$th column vector $\bx_\gamma$ are not
independent, the difference between the moments of the entries of the
standardized vector $\bx_\gamma$ and its unnormalized counterpart
$\widetilde\bx_\gamma$ is at least an order of magnitude smaller
than those
of~$\widetilde\bx_\gamma$. For instance, since $\bx_{i\gamma} =
O(N^{-1/2})$ for $1 \leq i \leq M$, for two independent ensembles of
covariance matrices $\widetilde{X}{}^{\bv}$ and $\widetilde{X}{}^{\bw}$
satisfying
(\ref{eqnXmat}) and~(\ref{eqnXmatexpbd}), we have the bound
%
%
\begin{equation}
\label{eqncovdiffbd} \E\bigl(\widetilde{ \bx}{}^{ \bv}_{i\gamma}
\bigr)^3-\E\bigl(\widetilde{\bx}{}^{
\bw}_{i\gamma
}
\bigr)^3=O\bigl(N^{-3/2}\bigr).
\end{equation}
On the other hand, if $\widetilde\bx_\gamma$ is the unnormalized
counterpart of $ \bx_\gamma$, as shown in Lem\-ma~\ref{lemEx},
%
%
\begin{equation}
\label{eqncorrdiffbd} \E(\widetilde{\bx}_{i\gamma})^3-\E({
\bx}_{i\gamma})^3 =O\bigl(N^{-5/2}\bigr).
\end{equation}
The above observation combined with a resolvent expansion---detailed in
Lemmas~\ref{2qmj},~\ref{lemYZ} and~\ref{lemEx}---gives
(\ref{ss440}).


\section{Proof of the main result}\label{sec4} \label{secpromainre}
In this section we will prove~(\ref{ss440}) in the case $X_0 = X^{ \bv
} = X$ and $X_{N} = X^{\bw} = \widetilde X$. As discussed above, it implies
(\ref{maincomp}) in Theorem~\ref{GFCT}. Similarly, one can prove
(\ref{c1}) and~(\ref{c11}) in Theorems~\ref{GFCT} and~\ref{GFCT2},
which complete the proof of Theorem~\ref{thmmain}, the main result of
this paper.

It is easy to see that~(\ref{ss440}) is a direct consequence of the
following lemma.\vadjust{\goodbreak}

%
\begin{lemma}\label{premainlemma}
Let $X $ be a $M \times N$ random matrix whose columns satisfy the
large deviation bounds~(\ref{eqnld1}),~(\ref{eqnld2})
and~(\ref{eqnld3}), for any $A_i \in\mathbb{C}$ and $B_{ij} \in
\mathbb{C}$ and for any $\zeta> 0$. The columns of $X$ are assumed to
be mutually independent. Furthermore, assume that the first column is
given by
%
%
\begin{equation}
\label{eqndeflemxi} X_{i1}= \frac{\widetilde x_{i1} }{\|\widetilde\bx_1
\|_2},\qquad 1\leq i \leq M,
%
\end{equation}
where $\widetilde x_{i1}$ are i.i.d. random variables with mean zero and
variance $M^{-1}$ and have an exponentially decay in the tails as given by
(\ref{eqnXmatexpbd}).

Let $\widetilde X$ be the random matrix whose entries have the same
distribution as~$X$ except for the first column, and the first column
of $\widetilde X$ is given by
\[
\widetilde X_{i1}=\widetilde x_{i1},
\]
where $\widetilde x_{i1}$ are as in~(\ref{eqndeflemxi}). The columns
of $\widetilde
X$ are also assumed to be mutually independent. Let $m,\widetilde m$ denote
the empirical Stieltjes transforms of $X^\dagger X$, $\widetilde
X^\dagger\widetilde X$.

Then for any function $F$ satisfying~(\ref{gflowder}), there exists
$\delta>0$, ${\varepsilon}_0>0$ depending only on $C_1$ such that for
any ${\varepsilon}<{\varepsilon}_0$
and for any real number $E$ satisfying
%
%
\begin{equation}
\label{Eeta} |E-\lambda_+|\leq N^{-2/3+\eps},\qquad {\eta_0} =
N^{-2/3-{\varepsilon}},
\end{equation}
we have
%
%
\begin{equation}
\label{maincomp3} \bigl|\E F \bigl( N \eta_0 \Im m (z) \bigr) - \E F
\bigl( N \eta_0 \Im\widetilde m (z) \bigr) \bigr| \leq O_{\varepsilon}
\bigl(N^{-1-\delta}\bigr),\qquad z=E+i\eta_0.\hspace*{-28pt}
\end{equation}
\end{lemma}
Note: In this lemma $X$ and $\widetilde X$ are neither pure correlation nor
pure covariance matrices, but their respective first columns are
distributed according to the standardized data matrix and raw data matrix.

%
\begin{remark}\label{remetasize} Under condition~(\ref{Eeta}) (see
\cite{NYcov}), we have the bound
%
%
\begin{equation}\quad
\label{immw}\quad C^{-1} \leq\bigl|m_W(z)\bigr| \leq C,\qquad \im
m_W (z) =O_{\varepsilon
}\bigl(N^{-1/3}\bigr),\qquad z = E + i
\eta_0.
\end{equation}
\end{remark}

First we collect some properties on submatrices of a \textit{generic}
$M\times N$ matrix~$Q$ which can be proved using standard results from
linear algebra.
Let~$Q^{(1)}$ be the $M \times(N-1)$ matrix obtained by removing the
first column of~$Q$. Define
%
%
\begin{equation}
\label{Gg110} G_Q^{(1)}= \bigl(\bigl(Q^{(1)}
\bigr)^\dagger\bigl(Q^{(1)}\bigr)-z \bigr)^{-1},\qquad
\mG_Q^{(1)}= \bigl(\bigl(Q^{(1)}\bigr)
\bigl(Q^{(1)}\bigr)^\dagger-z \bigr)^{-1}.
\end{equation}
Then\vspace*{1pt} by definition, $G_Q^{(1)}$ is a $(N-1)\times(N-1)$ matrix, $\mG
_Q^{(1)}$ is a $M\times M$ matrix and we have the identity
%
%
\begin{equation}
\label{propG1} \tr G_Q^{(1)}(z)-\tr\mG_Q^{(1)}(z)=
\frac{M-N+1}{ z}.
\end{equation}
Using the Cauchy interlacing theorem (see Equation (8.5) of
\cite{ESYY}), it can be shown that
%
%
\begin{equation}
\label{propG2} \tr G_Q^{(1)}(z)-\tr G_Q(z)
=O\bigl(\eta^{-1}\bigr),\qquad \eta=\im z.
\end{equation}

\begin{pf*}{Proof of Lemma~\ref{premainlemma}} First we note that
from Theorem 1.5 of~\cite{NYcov}, the conclusions of Theorem \ref
{lemstrMplawri} hold for both $X$ and $\widetilde X$.

Let $X^{(1)}$ be the $M \times(N-1)$ matrix obtained by removing the
first column of~$X$. Define
%
%
\begin{equation}\quad
\label{Gg11} G^{(1)}= \bigl(\bigl(X^{(1)}\bigr)^\dagger
\bigl(X^{(1)}\bigr)-z \bigr)^{-1},\qquad \mG^{(1)}= \bigl(
\bigl(X^{(1)}\bigr) \bigl(X^{(1)}\bigr)^\dagger-z
\bigr)^{-1}
\end{equation}
and as in~(\ref{dwbb}) set
%
%
\begin{equation}
\label{eqnmu} \mu= \eta_0 \im\tr G^{(1)} - \im
\frac{\eta_0}{z}.
\end{equation}
We will first verify that
%
%
\begin{eqnarray}
\label{ss44}\qquad
&&
\E F ( \eta_0 \im\tr G ) - \E F ( \mu) \nonumber\\
&&\qquad= \E
F^{(1) } ( \mu) (\im y_1+\im y_2+\im
y_3) +\E F^{(2) } ( \mu) \bigl(\tfrac12(\im
y_1)^2+\im y_1\im y_2\bigr)
\\
&&\qquad\quad{}+\E F^{(3) } ( \mu) \bigl(\tfrac16(\im y_1)^3
\bigr)+ O_{\varepsilon}\bigl(N^{-4/3}\bigr),
\nonumber
\end{eqnarray}
where $F^{(s)}$ denotes the $s$th derivative of $F$ and
$y_k$'s are defined as
%
%
\begin{equation}\label{eqndefyk}
y_k:=\eta_0 z \mw(-B)^{k-1} \bigl(
\bx_{1 },\bigl( \mG^{(1)}\bigr)^2
\bx_1 \bigr),
\end{equation}
where $ \bx_1 $ denotes the first column of $X$. Define the quantity
%
%
\begin{equation}
\label{215a} B:=-z m_W \biggl[ \bigl( \bx_1,
\mG^{(1)}(z) \bx_1 \bigr)- \biggl(\frac{-1}{z m_W(z)}-1
\biggr) \biggr].
\end{equation}
First, recall the following identity (see (6.23) of~\cite{NYcov}):
%
%
\begin{eqnarray}\quad
\label{26mm} \tr G -\tr G^{(1)} +z^{-1}
& = &
\bigl(G_{11}+z^{-1} \bigr)+\frac{ ( \bx_1, X^{(1)}G^{(1)}
G^{(1)}X^{(1)\dagger} \bx_1 ) }{- z - z ( \bx_1, \mG^{(1)}(z) \bx_1)}
\nonumber\\[-8pt]\\[-8pt]
&=& z G_{11} \bigl( \bx_{1 },\bigl( \mG^{(1)}
\bigr)^2 (z) \bx_1 \bigr).
\nonumber
\end{eqnarray}
Furthermore, as proved in Lemma 2.5 of~\cite{NYcov},
%
%
\begin{eqnarray}\label{eqnGii}
G_{11}(z) &=& \frac{ 1 }{- z - z ({\mathbf x}_1, \mG^{(1)}(z)
{\mathbf x}_1 ) }\quad \mbox{that is}\nonumber\\[-8pt]\\[-8pt]
\bigl( {\mathbf x}_1,
\mG^{(1)}(z) {\mathbf x}_1\bigr)&=&\frac{-1}{z G_{11}(z)}-1.\nonumber
\end{eqnarray}

From~(\ref{215a}) and~(\ref{eqnGii}) we obtain that
\[
B =-z m_W \biggl[ \biggl(\frac{-1}{z G_{11}(z)}-1 \biggr)- \biggl(
\frac{-1}{z m_W(z)}-1 \biggr) \biggr]=\frac{m_W-G_{11}}{G_{11}}.
\]
Fix $\zeta> 0$. From~(\ref{Lambdaofinal}), Remark~\ref{remetasize}
and the bound $|G_{11}| \leq|\mw| + O(1)$, it follows that for
$z = E + i \eta_0$,
%
%
\begin{equation}
\label{216a} | B| = \frac{ |\mw- G_{11}| }{|G_{11}|} \leq O_{\varepsilon}
\bigl(N^{-1/3}\bigr)\ll1
\end{equation}
with $\zeta$-high probability (see Definition~\ref{defhp}).
Therefore, with $\zeta$-high probability, we have the identity
%
%
\begin{equation}
\label{217a} G_{11}=\frac{m_{W}}{ B+1}=m_W\sum
_{k\geq0 }(-B)^k.
\end{equation}
Define $y$ to be the l.h.s. of~(\ref{26mm}) multiplied by $\eta_0$,
that is,
\[
y =\eta_0 \bigl( \tr G - \tr G^{(1)} + z^{-1}
\bigr),
\]
so that using~(\ref{26mm}) and~(\ref{217a}), we obtain
\[
y =\eta_0 z G_{11} \bigl( \bx_{1 },\bigl(
\mG^{(1)}\bigr)^2 \bx\bigr) = \sum
_{k=1}^\infty y_k.
\]
Since $ \bx_1$ satisfies~(\ref{eqnld1}),~(\ref{eqnld2}) and (\ref
{eqnld3}), and $\mG^{(1)}$ is independent of $\bx_1$, using Lemma
\ref{lemlargdev}, we infer that
for some $C_\zeta>0$
%
%
\begin{equation}
\label{hgjt} \bigl\llvert\bigl( \bx_{1 },\bigl( \mG^{(1)}
\bigr)^2 \bx_1 \bigr)\bigr\rrvert\leq\frac1M\tr\bigl(
\mG^{(1)}\bigr)^2+\frac{\varphi^{C_\zeta}}{M }\sqrt{ \tr\bigl|
\mG^{(1)}\bigr|^4}
\end{equation}
with $\zeta$-high probability. Using its definition, we bound $\tr(
\mG^{(1)})^2$ as
%
%
\begin{eqnarray}
\label{cG1} \bigl|\tr\bigl( \mG^{(1)}\bigr)^2\bigr|&\leq&\tr\bigl|
\mG^{(1)}\bigr|^2= \frac{\im\tr\mG^{(1)}}{ \eta_0}\nonumber\\[-8pt]\\[-8pt]
&=&O_{\varepsilon}
\bigl(N^{4/3}\bigr)+\frac{\im\tr G}{\eta_0} 
=O_{\varepsilon}
\bigl(N^{4/3}\bigr),\nonumber
\end{eqnarray}
where for the last two inequalities we have used~(\ref{propG1}),
(\ref{propG2}),~(\ref{Lambdafinal}) and~(\ref{immw}). Similarly, we
bound the last term of~(\ref{hgjt}) with
%
%
\begin{equation}
\label{cG2} \tr\bigl|\mG^{(1)}\bigr|^4\leq\eta_0^{-2}
\tr\bigl|\mG^{(1)}\bigr|^2\leq O_{\varepsilon}\bigl(N^{8/3}
\bigr)
\end{equation}
and obtain that
\[
\bigl\llvert\bigl( \bx_{1 },\bigl( \mG^{(1)}
\bigr)^2 \bx_1 \bigr)\bigr\rrvert\leq O_{\varepsilon}
\bigl(N^{1/3}\bigr).
\]
Equation~(\ref{216a}) and the fact $|z| + |m_{W}(z)|=O(1)$ yields that
%
%
\begin{equation}
\label{boundyk} |y_{k}|\leq O_{\varepsilon}\bigl(N^{-k/3 }
\bigr) \quad\mbox{and}\quad|y|\leq O_{\varepsilon}\bigl(N^{-1/3 }\bigr)
\end{equation}
holds with $\zeta$-high probability. Consequently, using (\ref
{gflowder}) and~(\ref{26mm}), we see that the expansion
%
%
\begin{equation}\qquad
\label{jkall} F ( \eta_0 \im\tr G ) - F (\mu) = \sum
_{k=1}^3 \frac{1 }{ k!} F^{(k)} \bigl(
N\eta_0 \Im\widetilde m^{(1)}(z) \bigr) (\Im
y)^k+O_{\varepsilon}\bigl(N^{-4/3}\bigr)\hspace*{-28pt}
\end{equation}
holds with $\zeta$-high probability. From the bounds on $y_k$'s
obtained above, equation
(\ref{ss44}) follows.

Now we estimate $\widetilde G$, which is defined as
\[
\widetilde G=\bigl(\widetilde X^\dagger\widetilde X-z
\bigr)^{-1}.
\]
Let $\widetilde{X}{}^{(1)}$ be the $M \times(N-1)$ matrix obtained by removing
the first column of~$\widetilde X$ and
$\widetilde\bx_1$ denote its first column.
Proceeding as in the previous calculations,
%
%
\begin{eqnarray}
\label{eqnlemmaGbytG}\qquad
&&
\E F ( \eta_0 \im\tr\widetilde G ) - \E F (
\mu) \nonumber\\
&&\qquad= \E F^{(1) } ( \mu) (\im\widetilde y_1+ \im
\widetilde y_2+\im\widetilde y_3 ) +\E F^{(2) } (
\mu) \bigl(\tfrac12(\im\widetilde y_1)^2+\im\widetilde
y_1\im\widetilde y_2 \bigr)
\\
&&\qquad\quad{} +\E F^{(3) } (\mu) \bigl(\tfrac16(\im\widetilde y_1)^3
\bigr)+O_{\varepsilon}\bigl(N^{-4/3}\bigr),
\nonumber
\end{eqnarray}
where
\begin{eqnarray*}
\widetilde y_k &=& \eta_0 z \mw(-\widetilde
B)^{k-1} \bigl(\widetilde\bx_{1 },\bigl( \mG^{(1)}
\bigr)^2 \widetilde\bx_1 \bigr),
\\
\widetilde B &=& -z m_W \biggl[ \bigl( \widetilde
\bx_1, \mG^{(1)}(z) \widetilde\bx_1 \bigr)-
\biggl(\frac{-1}{z m_W(z)}-1 \biggr) \biggr].
\end{eqnarray*}
Notice that $\mu$ appears in~(\ref{eqnlemmaGbytG}) because the
entries of $\widetilde X^{(1)}$ and $X^{(1)}$ are assumed to be identically
distributed.

Define the matrices
%
%
\begin{equation}
\label{defYZ} Y =\bigl(\mG^{(1)}\bigr)^2,\qquad
Z = \mG^{(1)}.
\end{equation}
The symmetric matrices $Y$ and $Z$ are independent of $
\bx_1$ and $\widetilde\bx_1$. Clearly, $YZ = ZY$.
Therefore, using the fact that $z$, $m_W\sim1$, we can write
\[
y_k=\eta_0 \sum_{0\leq n< k}
C_{k,n} ( \bx_1, Y \bx_1 ) (
\bx_1, Z \bx_1 )^n,
\]
where $C_{k,n} = O(1)$. Let $\cal Y= ( \bx_1, Y \bx_1
)$ and $\cal Z = ( \bx_1, Z \bx_1 )$. Then (\ref
{ss44}) can be written as
%
%
\begin{eqnarray}
\label{eqnkeylemfinest}
&&
\E F (\eta_0 \im\tr G ) - \E F ( \mu) \nonumber\\
&&\qquad=
\E F^{(1) } ( \mu)\im\biggl(\eta_0 \sum
_{0\leq n< k\leq3 } C_{k,n} \cal Y \cal Z^n \biggr)
\nonumber\\
&&\qquad\quad{}+ \E F^{(2) } ( \mu) \eta_0^2 \biggl(\frac12
\bigl(\im(C_{1,0} \cal Y ) \bigr)^2 +\im(C_{1,0}
\cal Y )\im(C_{2,0} \cal Y )\\
&&\qquad\quad\hspace*{128.5pt}{}+\im(C_{1,0} \cal Y )\im
(C_{2,1} \cal Y\cal Z ) \biggr)
\nonumber\\
&&\qquad\quad{}+\E F^{(3) } (\mu) \eta_0^3 \biggl(\frac16
\bigl(\im(C_{1,0} \cal Y ) \bigr)^3 \biggr)+O_{\varepsilon}
\bigl(N^{-4/3}\bigr).
\nonumber
\end{eqnarray}
Define $\widetilde{\cal Y}= ( \widetilde\bx_1, Y \widetilde
\bx_1 )$ and
$\widetilde{\cal Z} = ( \widetilde\bx_1, Z \widetilde\bx_1 )$. Using
(\ref{eqnlemmaGbytG}) and proceeding similarly as before, we obtain
that~(\ref{eqnkeylemfinest}) also holds for the case when $G$, $\cal
Y$ and $\cal Z$ are replaced with $\widetilde G$, $\widetilde{\cal Y}$
and $\widetilde
{\cal Z}$, respectively. The following is the key technical lemma of
this paper whose proof is deferred to the next section. 
%
%
\begin{lemma}\label{2qmj} Let $f\dvtx\R\to\R$ be a function satisfying
%
%
\begin{equation}
\label{gflowderf} \max_{x}\bigl|f (x)\bigr| \bigl(|x|+1 \bigr)^{-C } \leq C
\end{equation}
for some constant $C$. Let $\cal A$ be of the form
%
%
\begin{equation}
\label{byjc} \eta_0^a \prod
_{i=1}^a ( \bx, Y_i \bx) \prod
_{j=1}^b ( \bx, Z_j \bx),
\end{equation}
where $Y_i= Y$ or $Y^*$ and $Z_j= Z$ or $Z^*$ with $Y,Z$ as defined in
(\ref{defYZ}) and
$a, b$ are integers with $1\leq a\leq3, 1\leq a+b\leq3$. Then,
under the assumptions of Lem\-ma~\ref{premainlemma}, we have
%
%
\begin{equation}
\label{ldgg} \bigl|\E\bigl(f( \mu)\cal A \bigr)- \E\bigl( f( \mu
)\widetilde{\cal A}
\bigr) \bigr| \leq O_{\varepsilon}\bigl(N^{-7/6}\bigr),
\end{equation}
where $\widetilde{\cal A}$ is obtained by replacing $\bx$ with
$\widetilde\bx$ in
(\ref{byjc}).
\end{lemma}

Taking\vspace*{2pt} the difference of~(\ref{eqnkeylemfinest})
and the equation obtained by replacing~(\ref{eqnkeylemfinest}) with
$\widetilde G$, $\widetilde{\cal Y}$ and $\widetilde{\cal Z}$, we
deduce that the difference
\[
\E F (\eta_0 \im\tr G ) - \E F ( \eta_0 \im\tr
\widetilde G )
\]
can be\vspace*{2pt} approximated by the sum of $O(1)$ number of terms
of the form $\E (f( \mu)\cal A )- \E( f( \mu)\widetilde{\cal A} )$,
where $\cal A$ is as in~(\ref{byjc}) and $f$ is equal to
$F^{(1)}$, $F^{(2)}$ and $F^{(3)}$. 
Therefore, by applying Lemma~\ref{2qmj}, we conclude that Lemma \ref
{premainlemma} holds with any $\delta< 1/6$ and the proof is finished.
\end{pf*}

Finally, we are ready to give the proof of the main result of this paper:
\begin{pf*}{Proof of Theorem~\ref{thmmain}} By the Green function
comparison theorem discussed in Section~\ref{secproofsketch}, it only
remains to prove that
Theorems~\ref{GFCT} and~\ref{GFCT2} hold for the case
\[
X^\bv=X,\qquad X^\bw=\widetilde X.
\]
%
For simplicity, we will only prove~(\ref{maincomp}) of Theorem \ref
{GFCT}; the rest can be proved using almost identical arguments.

For $1\leq\gamma\leq N$, let $X_\gamma$ denote the random matrix
whose $j$th column is the same as that of $X^{\bv}$ if $j\geq\gamma$
and that of $X^{\bw}$
otherwise; in particular, $X_0 = X^{ \bv}$ and $X_{N} = X^{\bw}$. As
before, we define
\[
m_\gamma(z) =\frac{1}{N}\tr G_\gamma(z),\qquad
G_\gamma(z)=\bigl(X_\gamma^\dagger X_\gamma-z
\bigr)^{-1}.
\]
We have the telescoping sum,
%
%
\begin{eqnarray}
\label{eqntelesum1}
&&
\E^\bw F \bigl( N \eta_0 \Im
m^{\bw} (z) \bigr) - \E^\bv F \bigl( N \eta_0
\Im m^{\bv} (z) \bigr) \nonumber\\[-8pt]\\[-8pt]
&&\qquad= \sum_{\gamma= 1}^N
\E F \bigl( N \eta_0 \Im m_{\gamma} (z) \bigr) - \E F \bigl( N
\eta_0 \Im m_{\gamma-1} (z) \bigr).\nonumber
\end{eqnarray}
Applying Lemma~\ref{premainlemma} on $X_\gamma$ and $X_{\gamma-1}$
gives the estimate
%
%
\begin{equation}
\label{maincomp22}
\bigl|\E F \bigl( N \eta_0 \Im m_{\gamma} (z)
\bigr) - \E F \bigl( N \eta_0 \Im m_{\gamma-1} (z) \bigr) \bigr| \leq
O_{\varepsilon}\bigl(N^{-1-\delta}\bigr)
\end{equation}
for some $\delta> 0$.
Now~(\ref{maincomp}) follows from~(\ref{eqntelesum1}) and (\ref
{maincomp22}) and the
proof is finished.
\end{pf*}

%
\section{Moment computations}\label{sec5} \label{secmomcomp} In this
section we
prove Lemma~\ref{2qmj}.
For notational convenience, let us denote $\bx= \bx_1, \widetilde\bx
= \widetilde
\bx_1$.
We will also write
\[
\bx(k) = x_{k1},\qquad \widetilde\bx(k) = \widetilde x_{k1},\qquad
1 \leq k \leq M.
\]
Recall $\mu$ from~(\ref{eqnmu}). For the rest of this section, $a,b$
will denote
two integers with
\[
1\leq a\leq3,\qquad 1\leq a+b\leq3.
\]
Before stating the key results of this section, let us first give some
definitions.

%
\begin{definition}[{[$\I(A,\bk)$]}]
For any partition $A$ of the set $\{1,2,\ldots, 2a+2b\}$, and a vector
$\bk= \{k_1,k_2,\ldots, k_{2a+ 2b}\}, k_i \in\{1,2,\ldots,M\}$,
define the binary function $\I(A,\bk)$ as follows.
The function $\I(A,\bk)$ is equal to 1 if (1) for any $i,j$ in the
same block of $A$ we have $k_i=k_j$, (2) if $i,j$ are in different
blocks of $A$, we have $k_i\neq k_j$; otherwise $\I(A,\bk)=0$.
\end{definition}

%
\begin{example} If
%
%
\begin{equation}
\label{z37}
A=\bigl\{\{1\}, \{2,4\}, \{3,5,6\}\bigr\}
\end{equation}
and $a+b = 3$, then
\[
\I(A,\bk) ={\mathbf1}( k_2=k_4) {\mathbf1}( k_3=k_5=k_6){
\mathbf1}(k_1\neq k_2){\mathbf1}(k_2\neq
k_3){\mathbf1}(k_1\neq k_3).
\]
\end{example}
%
%
\begin{definition}[{[$\mathcal{N}(A,1), \mathcal{N}(A,2)$ and
${\mathbf I}_{(A,3)}$]}] Given a
partition $A$ of the set $\{1,2,\ldots, 2a+2b\}$, let $\mathcal
{N}(A,1)$ be the
number of the blocks in $A$ that contain only one element of the set $\{
1,2 ,\ldots, 2a+2b\}$. Let $\mathcal{N}(A,2)$ be the number of the
blocks in $A$ of the form $\{k_{2i-1}, k_{2i}\}$ with $i > a$. Note
that $\mathcal{N}(A,2)$ depends on $a$ and $b$ in addition to $A$. Let
${\mathbf I}_{(A,3)}$ be equal to one if and only if $a+b=3$ and $A$ is
composed of
$2$ blocks with three elements in each block.
\end{definition}


The proof of Lemma~\ref{2qmj} relies on Lemmas~\ref{lemYZ} and \ref
{lemEx} stated below and proved at the end of this section.\vadjust{\goodbreak}
%
%
\begin{lemma}\label{lemYZ} Recall the matrices $Y,Z$ from (\ref
{defYZ}). Then for any ${\varepsilon}>0$ the following estimate
\begin{eqnarray*}
&&\sum_{k_1,k_2,\ldots, k_{2a+2b} =1}^M \I(A,\bk)
\eta^a_0 (Y_{k_1k_2}\cdots Y_{k_{2a-1}k_{2a}} )
(Z_{k_{2a+1}k_{2a+2}}\cdots Z_{k_{2a+2b-1}k_{2a+2b}} )
\\
&&\qquad=O_{\varepsilon} \bigl(\bigl( N^{2 /3}\bigr)^{a+b}
\bigl(N^{1/2}\bigr)^{\mathcal
{N}(A,1)+{\mathbf I}_{(A,3)}
}\bigl(N^{1/3}
\bigr)^{\mathcal{N}(A,2)} \bigr)
\end{eqnarray*}
holds with $\zeta$-high probability for any fixed $\zeta>0$. The
result also holds if
any of the $Y,Z$ are replaced by their complex conjugates
$Y^*,Z^*$, respectively.
\end{lemma}

%
\begin{lemma}\label{lemEx}
Let $\widetilde y_i$ be i.i.d. random variables such that
\[
\E\widetilde y_i=0,\qquad \E(\widetilde y_i
)^2 = M^{-1},\qquad 1 \leq i \leq M,
\]
and have a subexponential decay as in~(\ref{eqnXmatexpbd}).
Let $A$ be a partition of the set $\{1,2 ,\ldots, 2a+2b\}$ and let
\[
y_i:=\frac{\widetilde y_i }{(\sum_j\widetilde y_j^{ 2})^{ 1/2}}.
\]
Then for any vector $\bk= (k_1,k_2,\ldots,k_{2a+2b})$ and for any
${\varepsilon}> 0$, we have
%
%
\begin{eqnarray}
\label{x11}
&&
\E\Biggl( \I(A,\bk)\prod_{i=1}^{2a+2b}
y_{k_i} \Biggr) - \E\Biggl( \I(A,\bk) \prod
_{i=1}^{2a+2b} \widetilde y_{k_i}
\Biggr)\\
&&\qquad=O_{\varepsilon
}\bigl(N^{-(a+b)-\max\{\mathcal{N}(A,1),1\}}\bigr).\nonumber
\end{eqnarray}
\end{lemma}

With the above two lemmas in hand, we are now ready to give the proof of
Lemma~\ref{2qmj}.
\begin{pf*}{Proof of Lemma~\ref{2qmj}}
We will only prove the case when
%
%
\begin{equation}
\label{YiZi} Y_i=Y,\qquad Z_i=Z
\end{equation}
for all $i$
and, thus,
\[
\cal A= \eta_0^a ( \bx, Y \bx)^a ( \bx, Z
\bx)^b.
\]
The other cases can be proved similarly. First, let us write (\ref
{byjc}) as
\begin{eqnarray*}
&&\eta_0^a ( \bx, Y \bx)^a( \bx, Z \bx
)^b
\\
&&\qquad=\sum_{A}\sum_{k_1,k_2,\ldots, k_{2a+2b} =1}^M
\eta_0^a \I(A,\bk) 
\prod_{i=1}^{2a+2b}
\bx(k_i) (Y_{k_1k_2}\cdots Y_{k_{2a-1}k_{2a}} )\\
&&\qquad\quad\hspace*{150.5pt}{}\times(Z_{k_{2a+1}k_{2a+2}}\cdots Z_{k_{2a+2b-1}k_{2a +2b}} ),
\end{eqnarray*}
where the summation index $A$ ranges over all the partitions of the set
$\{1,2,\ldots, 2a+2b\}$. Taking expectations, and using the fact that
$\bx$ is independent of $Y$, $Z$ and~$\mu$, leads to
%
%
\begin{eqnarray}
\label{19bb} \E f( \mu) \cal A &=& \sum_A \sum
_{k_1,k_2,\ldots, k_{2a+2b}=1}^M \E\Biggl(\eta_0^a
\I(A,\bk) \nonumber\\[-2pt]
&&\qquad\quad\hspace*{60pt}{}\times\prod_{i=1}^{2a+2b}
\bx(k_i) (Y_{k_1k_2}\cdots Y_{k_{2a-1}k_{2a}} )\nonumber\\[-2pt]
&&\hspace*{128.5pt}{}\times
(Z_{k_{2a+1}k_{2a+2}}\cdots Z_{k_{2a+2b-1}k_{2a+2b}} ) \Biggr)
\nonumber\\[-9pt]\\[-9pt]
&=& \sum_A \Biggl(\E\I(A,\bk) \prod
_{i=1}^{2a+2b} \bx({k_i}) \Biggr)
\nonumber\\[-2pt]
&&\hspace*{14.5pt}{}\times \Biggl( \E f( \mu)\sum_{k_1,k_2,\ldots, k_{2a+2b}=1}^M \I(A,
\bk) \eta_0^a (Y_{k_1k_2}\cdots Y_{k_{2a-1}k_{2a}}
)\nonumber\\[-2pt]
&&\hspace*{132.6pt}{} \times(Z_{k_{2a+1}k_{2a+2}}\cdots Z_{k_{2a+2b-1}k_{2a+2b}} )
\Biggr),\nonumber
\end{eqnarray}
where the last inequality follows from the fact that $ (\E \I
(A,\bk) \prod_{i=1}^{2a+2b} \bx({k_i}) )$ is independent of $Y,Z$.
Combining~(\ref{19bb}), Lemmas~\ref{lemYZ} and~\ref{lemEx},
we deduce that
%
%
\fontsize{10pt}{\baselineskip}\selectfont
\makeatletter
\def\tagform@#1{\normalsize\maketag@@@{(\ignorespaces#1\unskip\@@italiccorr)}}
\makeatother\begin{eqnarray}
\label{eqnkeylemcont}
&&
\bigl|\E\bigl(f( \mu)\cal A \bigr)- \E\bigl( f( \mu)
\widetilde{\cal A} \bigr) \bigr| \hspace*{-15pt}\nonumber\\
&&\qquad\leq \sum_{ A}O_{\varepsilon}
\bigl(\bigl( N^{-1 /3}\bigr)^{a+b} \bigl(N^{1/2}
\bigr)^{\mathcal{N}(A,1)+{\mathbf I}_{(A,3)}}\bigl(N^{ -1}\bigr)^{\max\{
\mathcal{N}(A,1),1\}}
\bigl(N^{1/3}\bigr)^{\mathcal
{N}(A,2) } \bigr)
\hspace*{-15pt}\\
&&\qquad\leq \sum_{ A}O_{\varepsilon} \bigl(\bigl(
N^{-a /3}\bigr) \bigl(N^{1/2}\bigr)^{\mathcal{N}(A,1)
+{\mathbf I}_{(A,3)}}
\bigl(N^{ -1}\bigr)^{\max\{\mathcal{N}(A,1),1\}
}\bigl(N^{1/3}
\bigr)^{\mathcal{N}(A,2) -b} \bigr).\hspace*{-15pt}
\nonumber
\end{eqnarray}
\normalsize
Now we claim that the terms in the r.h.s. of~(\ref{eqnkeylemcont})
are bounded by $O_{\varepsilon}(N^{-7/6})$. Indeed, note that
$\mathcal{N}(A,1)
>0$ implies ${\mathbf I}_{(A,3)}= 0$. Therefore, the worse case
scenario is the case
in which
\[
a=1,\qquad b=\mathcal{N}(A,2) \quad\mbox{and}\quad\mathcal{N}(A,1)=1,
\]
since by definition we have $\mathcal{N}(A,2) \leq b$.
But it is easy to see the above scenario cannot occur, since if the
first two conditions hold, then it follows that
$\mathcal{N}(A,1)=0$ or $2$. Thus, we have finished the proof of Lemma
\ref{2qmj}.
\end{pf*}

\begin{pf*}{Proof of Lemma~\ref{lemYZ}}
Note that all of the bounds
in this lemma hold with $\zeta$-high probability, not in expectation.
For simplicity, we
will subsume this in the notation.\vadjust{\goodbreak}

First let us prove a slightly different result. Define the binary function
$\widetilde\I(A,\bk)$ [similar to $\I(A,\bk)$] as follows.
$\widetilde\I(A,\bk
)$ is equal to 1 in the following scenarios: (1) for any $i,j$ in the
same block of $A$ we have $k_i=k_j$, (2) if $i,j$ are in different
blocks of $A$, we have $k_i\neq k_j$ except that if one of the indices
$i,j$ is in the block of $A$ which contains exactly two elements, then
$k_i$ is allowed to be equal to $k_j$. In all other instances
$\widetilde\I
(A,\bk)=0$. For instance, in the previous example~(\ref{z37}), we have
\[
\widetilde\I(A,\bk) ={\mathbf1}( k_2=k_4) {\mathbf1}(
k_3=k_5=k_6) {\mathbf1}(k_1\neq
k_3).
\]

We first claim that
%
%
\fontsize{11pt}{\baselineskip}\selectfont
\makeatletter
\def\tagform@#1{\normalsize\maketag@@@{(\ignorespaces#1\unskip\@@italiccorr)}}
\makeatother\begin{eqnarray}
\label{511qz0}
&&
\sum_{k_1,k_2,\ldots, k_{2a+2b}=1}^M \widetilde
\I(A,\bk) \eta_0^a (Y_{k_1k_2}\cdots
Y_{k_{2a-1}k_{2a}} ) (Z_{k_{2a+1}k_{2a+2}}\cdots
Z_{k_{2a+2b-1}k_{2a+2b}} )\hspace*{-22pt}
\nonumber\\[-6pt]\\[-10pt]
&&\qquad=O_{\varepsilon} \bigl(\bigl( N^{2 /3}\bigr)^{a+b}
\bigl(N^{1/2}\bigr)^{\mathcal
{N}(A,1)+ {\mathbf I}_{(A,3)}
}\bigl(N^{1/3}
\bigr)^{\mathcal{N}(A,2)} \bigr).
\nonumber
\end{eqnarray}
\normalsize
%
Let us first prove~(\ref{511qz0}) when ${\mathbf I}_{(A,3)}= 0$.
Define the functions
\[
g_1(m):= \tr\bigl|Z^m\bigr|,\qquad g_2(m):= \sqrt{
\bigl({\tr}|Z|^{2m}\bigr)},\qquad 1 \leq m \leq2a +b.
\]
%
We will show that the
%
%
\begin{equation}
\label{eqnlhswoex}
\mbox{l.h.s. of~(\ref{511qz0})} \leq
O_{\varepsilon} \biggl(\eta_0^a\bigl(N^{1/2}
\bigr)^{\mathcal{N}(A,1)}\prod_i g_{\al
_i}(m_i)
\biggr),
\end{equation}
where $ \alpha_{i} \in\{1,2\}$ and $ m_i \leq2a +b$.

To this end, we will use the following {2--1--3 rule}: 
%
\begin{itemize}
\item{2:} If the index $i$ appears in a block of $A$ which
contains exactly two elements, first sum up over the index $k_i$. Then
estimate the remaining terms with absolute sum. For example, let $A=\{\{
1\}, \{2,3\},\{4\}\}$. Recall that $Y=Z^2$,
\[
\biggl| \sum_{k_1,k_2,k_3} \widetilde\I(A,\bk)
Y_{k_1k_2}Z_{k_2k_4}\biggr|=\biggl|\sum_{k\neq l}
(YZ)_{kl}\biggr|\leq\sum_{kl}
\bigl|(YZ)_{kl}\bigr|=\sum_{kl} \bigl|\bigl(
Z^3\bigr)_{kl}\bigr|.
\]
\item{1:} Next do the summation over the index $k_i$ if $i$
appears in the block of $A$ which contains only one element as follows:
\[
\sum_{ l }\bigl|\bigl(Z^m
\bigr)_{kl} \bigr|\leq CN^{1/2} \sqrt{\bigl(|Z|^{2m}
\bigr)_{kk}},\qquad \sum_{ kl }\bigl|
\bigl(Z^m\bigr)_{kl} \bigr|\leq CN \sqrt{{\tr}|Z|^{2m}
}.
\]
In the above inequalities, we have used the Cauchy--Schwarz and the
fact that $Z$ is a symmetric matrix. Note that each summation of the
above kind brings an extra $N^{1/2}$ factor.
\item{3:} Finally, sum up over the other indices. After the first
two steps,~(\ref{511qz0}) will be reduced to the product of following
terms: 
%
\[
\bigl(N^{1/2}\bigr)^{\mathcal{N}(A,1)},\qquad \bigl|\tr Z^r\bigr|,\qquad %
\sqrt{{\tr}|Z|^{2r} },\qquad r \leq2a+ b,\vadjust{\goodbreak} 
\]
and terms of the form
%
%
\begin{equation}
\label{36sg} \sum_{k} \prod
_{i=1}^m\bigl|\bigl(Z^{m_i}
\bigr)_{kk}\bigr|\prod_{j=1}^n\sqrt
{\bigl(|Z|^{2n_j}\bigr)_{kk}},\qquad 2\leq m+n.
\end{equation}
%
If $m+n=2$, then using the Cauchy--Schwarz inequality,~(\ref{36sg})
can be estimated as
%
%
\begin{equation}\quad
\label{eqnmpne2o} \prod_{i=1}^m \prod
_{j=1}^n \sum
_{k} \bigl|\bigl(Z^{m_i}\bigr)_{kk}\bigr| \sqrt{
\bigl(|Z|^{2n_j}\bigr)_{kk}} \leq\prod
_{i=1}^m \prod_{j=1}^n
\sqrt{{\tr}|Z|^{2m_i} }\sqrt{{\tr}|Z|^{2n_j} }.
\end{equation}
For $m+n > 2$, we bound $m+n-2$ of them [$|(Z^{m_i})_{kk}|$ or $\sqrt
{(|Z|^{2n_j})_{kk}}$] by the maximum as follows:
\begin{eqnarray*}
\bigl|\bigl(Z^{m_i}\bigr)_{kk}\bigr|&\leq&\max_k \bigl|
\bigl(Z^{m_i}\bigr)_{kk}\bigr|\leq\sqrt{{\tr}|Z|^{2m_i}},\\
\sqrt{\bigl(|Z|^{2n_j}\bigr)_{kk}}&\leq&\max_k\sqrt{
\bigl(|Z|^{2n_j}\bigr)_{kk}}\leq\sqrt{{\tr}|Z|^{2n_j}},
\end{eqnarray*}
to reduce to the case of $m+n=2$ and use the bound~(\ref{eqnmpne2o}).
\end{itemize}
Let us give an example in the case $a=1$, $b=2$ and $A=\{\{1\}, \{2, 3\},
\{4, 5,6\}\}$. Then the term~(\ref{511qz0}) in this case reduces to
\[
\sum_{k_1k_2k_4}\eta_0 Y_{k_1k_2}Z_{k_2k_4}Z_{k_4k_4}
\leq\sum_{k_1 k_4}\eta_0 \bigl\llvert
\bigl(Z^3\bigr)_{k_1k_4}\bigr\rrvert\llvert Z_{k_4k_4}
\rrvert,
\]
where the above inequality is obtained by applying rule
{2}. Next, applying
rule {1} yields
\[
\leq\sum_{ k_4}\eta_0 N^{1/2}
\sqrt{\bigl(|Z|^6\bigr)_{k_4k_4}} \llvert Z_{k_4k_4}
\rrvert
\]
and, finally, applying rule {3} leads to the bound
\[
\leq\sum_{ k_4}\eta_0 N^{1/2}
\sqrt{{\tr}|Z|^6} \sqrt{{\tr}|Z|^2}.
\]
%
Using this {2--1--3} rule described above, we obtain~(\ref{eqnlhswoex}).
By the definition of the {2--1--3} rule, it is easy to see that
%
%
\begin{equation}
\label{515sh} \sum_i m_i= 2a+b.
\end{equation}
Recall $\eta_0= N^{-2/3-{\varepsilon}}$. Using~(\ref{cG1}) and
(\ref{cG2}), we deduce that
if $\al_i m_i\neq1$, then
\[
g_{\al_i}(m_i)\leq O_{\varepsilon}\bigl(N^{2m_i/3}
\bigr).
\]
For $\al_i m_i = 1$, using~(\ref{propG1}),~(\ref{propG2}), (\ref
{Lambdafinal}) and $m_W=O(1)$, we see that $g_{1}(1)=O_{\varepsilon
}(N)$. Thus,
%
%
\begin{equation}
\label{516sh} g_{\al_i}(m_i)\leq O_{\varepsilon}
\bigl(N^{2m_i/3}\bigr) \bigl(N^{1/3}\bigr)^{{\mathbf1}( \al
_i m_i=1) }.
\end{equation}
Combining equations~(\ref{eqnlhswoex})--(\ref{516sh}), we have the
%
%
\begin{eqnarray}
\label{xjfs}
&&\mbox{l.h.s. of equation~(\ref{511qz0})}\nonumber\\[-8pt]\\[-8pt]
&&\qquad =
O_{\varepsilon}\bigl(N^{1/2}\bigr)^{\mathcal
{N}(A,1)}N^{2a/3+2b/3}
\bigl(N^{1/3}\bigr)^{\#\{
i\dvtx a_im_i=1\}}.\nonumber
\end{eqnarray}
Now notice that by the definition, the term $g_1(1)$ in (\ref
{eqnlhswoex}) can only be created during the first step of the
{2--1--3} rule, that is, the {2} rule, and, therefore, we
deduce that
\[
\mathcal{N}(A,2) =\#\{i\dvtx\alpha_im_i=1\},
\]
which completes the proof of the claim made in~(\ref{511qz0}) for the
case ${\mathbf I}_{(A,3)}= 0$.

Now consider\vspace*{2pt} the case ${\mathbf I}_{(A,3)}= 1$. 
Using the fact that $Y,Z$ are symmetric matrices and the relation $Y = Z^2$,
we deduce that the term
\[
\sum\widetilde\I(A,\bk) (Y_{k_1k_2}\cdots Y_{k_{2a-1}k_{2a}} )
(Z_{k_{2a+1}k_{2a+2}}\cdots Z_{k_{2a+2b-1}k_{2a+2b}} )
\]
reduces to one of the following situations:
%
%
\begin{eqnarray}
\label{esmallfake}
&&
\sum_{k_1,k_2,\ldots, k_{2a+2b}}\widetilde\I(A,\bk)
(Y_{k_1k_2}\cdots Y_{k_{2a-1}k_{2a}} ) (Z_{k_{2a+1}k_{2a+2}}\cdots
Z_{k_{2a+2b-1}k_{2a+2b}} )\nonumber\hspace*{-35pt}\\[-8pt]\\[-8pt]
&&\qquad =
\cases{ \displaystyle \sum_{k_1,k_2 = 1}^M
Z^{\b_1}_{k_1k_1}Z^{\b_2}_{k_1k_2}Z^{\b
_3}_{k_2k_2},
\cr
\displaystyle \sum_{k_1,k_2 =
1}^MZ^{\b_1}_{k_1k_2}Z^{\b_2}_{k_1k_2}Z^{\b_3}_{k_1k_2},}\nonumber\hspace*{-35pt}
\end{eqnarray}
for $\b_i \in\{1,2\}, i \in\{1,2,3\}$.
We bound the first scenario above as
%
%
\begin{eqnarray}
&&\biggl\vert\sum_{k_1k_2} \bigl(Z^{\b_1}
\bigr)_{k_1k_1}\bigl(Z^{\b_2}\bigr)_{k_1k_2}
\bigl(Z^{\b
_3}\bigr)_{k_2k_2} \biggr\vert
\nonumber\\
&&\qquad\leq \sum_{k_1k_2} \bigl\vert\bigl(Z^{\b_1}
\bigr)_{k_1k_1}\bigl(Z^{\b_2}\bigr)_{k_1k_2} \bigr\vert
\max_k\bigl\llvert\bigl(Z^{\b_3}\bigr)_{kk}\bigr
\rrvert\\
&&\qquad\leq\sum_{k_1k_2} \bigl\vert\bigl(Z^{\b_1}
\bigr)_{k_1k_1}\bigl(Z^{\b_2}\bigr)_{k_1k_2} \bigr\vert\sqrt{
{\tr}|Z|^{2\b_3} }.
\nonumber
\end{eqnarray}
Using rule {1} and rule {3} above yields
\[
\sum_{k_1k_2} \bigl\vert\bigl(Z^{\b_1}
\bigr)_{k_1k_1}\bigl(Z^{\b_2}\bigr)_{k_1k_2} \bigr\vert\leq
CN^{1/2} \sqrt{{\tr}|Z|^{2\b_1} } \sqrt{{\tr}|Z|^{2\b_2} }
\]
and, thus,
\begin{eqnarray*}
\eta^0_a\sum_{k_1,k_2 = 1}^M
\bigl\vert Z^{\b_1}_{k_1k_1}Z^{\b
_2}_{k_1k_2}Z^{\b_3}_{k_2k_2}
\bigr\vert
&\leq& C \eta^0_a N^{1/2} \sqrt{
{\tr}|Z|^{2\b_1} } \sqrt{{\tr}|Z|^{2\b_2}}\sqrt{{\tr}|Z|^{2\b
_3}}
\\
&=& O_{\varepsilon}\bigl(N^{-2a/3 +{1/2}}\bigr) O_{\varepsilon}
\bigl(N^{2/3(\b_1 +
\b_2+\b_3)}\bigr) \\
&=& O_{\varepsilon}\bigl(N^{2/3(a+b) +{1/2}}\bigr),
\end{eqnarray*}
where in the last inequality we have used the fact that $\sum_i \b_i
= 2a+b$.
For the second case in~(\ref{esmallfake}), first we note
\[
\max_{kl} \bigl|\bigl(Z^{\b}\bigr)_{kl }\bigr|\leq\sqrt{
{\tr}|Z|^{2\b} }.
\]
Now using the Cauchy--Schwarz inequality,
\[
\sum_{k_1,k_2}\bigl\vert\bigl(Z^{\b_1}
\bigr)_{k_1k_2}\bigl(Z^{\b_2}\bigr)_{k_1k_2}
\bigl(Z^{\b
_3}\bigr)_{k_1k_2} \bigr\vert\leq\sqrt{
{\tr}|Z|^{2\b_1 }} \sqrt{{\tr}|Z|^{2\b_2 }} \sqrt{{\tr}|Z|^{2\b
_3 }}
\]
and, thus,
\begin{eqnarray*}
\eta^0_a\sum_{k_1,k_2 = 1}^M
\bigl\vert Z^{\b_1}_{k_1k_2}Z^{\b
_2}_{k_1k_2}Z^{\b_3}_{k_1k_2}
\bigr\vert
&\leq& C \eta^0_a \sqrt{{\tr}|Z|^{2\b_1} }
\sqrt{{\tr}|Z|^{2\b_2}}\sqrt{{\tr}|Z|^{2\b
_3}}
\\
&=& O_{\varepsilon}\bigl(N^{-2a/3}\bigr) O_{\varepsilon}
\bigl(N^{2/3(\b_1 + \b_2+\b
_3)}\bigr) \\
&=& O_{\varepsilon}\bigl(N^{2/3(a+b)}\bigr).
\end{eqnarray*}

Summarizing the above computations, and noticing that $\mathcal
{N}(A,1)=\break \mathcal
{N}(A,2)= 0$ when ${\mathbf I}_{(A,3)}= 1$, we obtain the bound
\begin{eqnarray*}
&&
\eta^a_0 \widetilde\I(A,\bk) \bigl\vert(Y_{k_1k_2}
\cdots Y_{k_{2a-1}k_{2a}} ) (Z_{k_{2a+1}k_{2a+2}}\cdots
Z_{k_{2a+2b-1}k_{2a+2b}} ) \bigr\vert\\
&&\qquad=
O_{\varepsilon}\bigl(N^{2/3(a+b)
+{1/2}}\bigr)
\\
&&\qquad= O_{\varepsilon} \bigl(\bigl( N^{2 /3}\bigr)^{a+b}
\bigl(N^{1/2}\bigr)^{{\mathbf
I}_{(A,3)}} \bigr),
\end{eqnarray*}
proving the claim~(\ref{511qz0}) when ${\mathbf I}_{(A,3)}= 1$.

Now we return to prove Lemma~\ref{lemYZ}. One can see that for any
partition~$A$ of the set $\{1,2,\ldots, 2a + 2b\}$ and a vector $\bk
$, the function $\I(A,\bk)$ can be written as linear combinations of
the functions~$\widetilde\I(A_i,\bk)$ for some partitions~$A_i$'s of the
set $\{1,2,\ldots, 2a + 2b \}$ such that
\[
\mathcal{N}(A_i,1) \leq\mathcal{N}(A,1),\qquad \mathcal{N}(A_i,2)
\leq\mathcal{N}(A,2),\qquad {\mathbf I}_{A_i,3}={\mathbf I}_{(A,3)}.
\]
For instance, for $A$ given in~(\ref{z37}),
\[
\widetilde\I(A,\bk) ={\mathbf1}( k_2=k_4) {\mathbf1}(
k_3=k_5=k_6) {\mathbf1}(k_1\neq
k_3),
\]
we have the identity
\[
\I(A,\bk)= \widetilde\I(A,\bk)- \widetilde\I(A_1,\bk)- \widetilde
\I(A_2,\bk),
\]
where $A_1=\{\{1\}, \{2,3,4,5,6\}\}$ and $A_2=\{\{1,2,4\}, \{ 3, 5,6\}\}
$. Now the lemma follows from~(\ref{511qz0}) and the proof is finished.
\end{pf*}
\begin{pf*}{Proof of Lemma~\ref{lemEx}}
For any $k_1, k_2, \ldots, k_m \in\{1,2,\ldots, M\}$ and \mbox{$m \in
\mathbb{N}$}, by definition we have
%
%
\begin{eqnarray}
\label{eqntayexp} \E\I(A,\bk) \prod_{i=1}^{m}
y_{k_i} &=&\E\I(A,\bk) \frac{\prod_{i=1}^{m} \widetilde y_{k_i} }{(\sum_j
\widetilde
y_j^{ 2})^{m/2}} \nonumber\\[-8pt]\\[-8pt]
&=&\E\I(A,\bk){\prod
_{i=1}^{m} \widetilde y_{k_i}} {
\Biggl[1-\sum_{j=1}^M \biggl(
\frac{1 }{ M}-\widetilde y_j^{ 2}\biggr)
\Biggr]^{-m/2}}.\nonumber
\end{eqnarray}
Using large deviation bounds, it is easy to see that for any
${\varepsilon}>0$
%
%
\begin{equation}
\label{343} \sum_{j=1}^M \biggl(
\frac{1}{ M}-\widetilde y_j^{ 2}
\biggr)=O_{\varepsilon
}\bigl(N^{-1/2}\bigr).
\end{equation}
Therefore, by the Taylor expansion,
%
%
\begin{eqnarray}
\label{eqnkeylemcombf}\qquad
&&
\E\I(A,\bk) \prod_{i=1}^{2a+2b}
y_{k_i} - \E\I(A,\bk) \prod_{i=1}^{2a+2b}
\widetilde y_{k_i} 
\nonumber\\[-8pt]\\[-8pt]
&&\qquad= \sum_{n=1}^\infty C_{n} \E
\Biggl[ \I(A,\bk) \Biggl(\prod_{i=1}^{2a+2b}
\widetilde y_{k_i} \Biggr) \Biggl(\sum_{r_1,r_2,\ldots,r_n = 1}^M
\prod_{j=1}^n\biggl(\frac{1}{ M}-
\widetilde y_{r_j}^{ 2}\biggr) \Biggr) \Biggr],
\nonumber
\end{eqnarray}
where $C_{n} = C_{a,b,n}$ is a combinatorial factor. Using
(\ref{343}), the r.h.s. of equation~(\ref{eqnkeylemcombf}) may be
expressed as
%
%
\begin{eqnarray}\quad
\label{djs} &=& \sum_{n=1}^{n_0}
C_n \E\Biggl[ \I(A,\bk) \Biggl(\prod_{i=1}^{2a+2b}
\widetilde y_{k_i} \Biggr) \Biggl(\sum_{r_1,r_2,\ldots,r_n = 1}^M
\prod_{j=1}^n\biggl(\frac{1 }{ M} -
\widetilde y_{r_j}^{ 2}\biggr) \Biggr) \Biggr]\nonumber\\[-8pt]\\[-8pt]
&&{}+O_{\varepsilon}\bigl(\bigl(N^{-1/2}\bigr)^{2a+2b+n_0
}\bigr)\nonumber
\end{eqnarray}
for some fixed $n_0 \in\mathbb{N}$ (say, $n_0 =20$). 

Since $n_0, a, b = O(1)$,
the combinatorial factors do not increase with $N$, that is,
$C_{n} = O(1)$, and, thus, we can bound
%
%
\begin{equation}
\label{ddj} 
\E\Biggl[ \I(A,\bk) \Biggl(\prod
_{i=1}^{2a+2b} \widetilde y_{k_i} \Biggr)
\Biggl( \prod_{j=1}^n\biggl(
\frac{1 }{ M} -\widetilde y_{r_j}^{ 2}\biggr) \Biggr)
\Biggr]
\end{equation}
as follows.
Notice that the number of distinct indices $k_i$ in~(\ref{ddj}) is
equal to the number
of blocks in the partition $A$. Thus, for a given set of values for the
indices $r_1,r_2,\ldots, r_n$, the term~(\ref{ddj}) is nonzero only
if at least $\mathcal{N}(A,1)$ of the indices $r_j$ belong to the set
$\{k_1,k_2,\ldots, k_{2a + 2b}\}$. The above observation also implies
that for~(\ref{ddj}) to be nonzero we must have
%
%
\begin{equation}
\label{conc} n\geq\mathcal{N}(A,1).
\end{equation}
Furthermore, the indices $r_j$ which do not belong to the set $\{
k_1,k_2,\ldots, k_{2a + 2b}\}$ must appear more than once since $\E
(1/M - y^2_{r_j}) = 0$. This crucial observation implies that, if the
term~(\ref{ddj}) is nonzero and
%
%
\begin{equation}
\label{eqnno2na0} \mathcal{N}(A,1)=0
\qquad\mbox{then } n\geq2.
\end{equation}
Therefore, the number of nonzero terms in the sum
%
%
\begin{equation}
\label{ddjsum} \sum_{r_1,r_2,\ldots,r_n = 1}^M \E\Biggl[
\I(A,\bk) \Biggl(\prod_{i=1}^{2a+2b} \widetilde
y_{k_i} \Biggr) \Biggl( \prod_{j=1}^n
\biggl(\frac{1 }{ M} -\widetilde y_{r_j}^{ 2}\biggr)
\Biggr) \Biggr]
\end{equation}
is $O((N^{1/2})^{n - \mathcal{N}(A,1)})$, and each of these
terms are of the size $O_{\varepsilon}(N^{-(a+b) -n})$, yielding
%
%
\begin{eqnarray}
\label{eqnmclemkeybd}
&&
\sum_{r_1,r_2,\ldots,r_n = 1}^M \E
\Biggl[ \I(A,\bk) \Biggl(\prod_{i=1}^{2a+2b}
\widetilde y_{k_i} \Biggr) \Biggl( \prod_{j=1}^n
\biggl(\frac{1 }{ M} -\widetilde y_{r_j}^{ 2}\biggr)
\Biggr) \Biggr] \nonumber\\[-8pt]\\[-8pt]
&&\qquad\leq O_{\varepsilon}\bigl(N^{-(a+b) -n/2- \mathcal
{N}(A,1)/2}\bigr).\nonumber
\end{eqnarray}
Combining~(\ref{eqnmclemkeybd}) with~(\ref{conc}) and the
observation made in~(\ref{eqnno2na0}),
we obtain that
%
%
\begin{eqnarray}
\label{hddj}
&&\sum_{r_1,r_2,\ldots,r_n = 1}^M \E\Biggl[
\I(A,\bk) \Biggl(\prod_{i=1}^{2a+2b} \widetilde
y_{k_i} \Biggr) \Biggl( \prod_{j=1}^n
\biggl(\frac{1 }{ M} -\widetilde y_{r_j}^{ 2}\biggr)
\Biggr) \Biggr] \nonumber\\[-8pt]\\[-8pt]
&&\qquad\leq O_{\varepsilon}\bigl(N^{-(a+b)-\max\{\mathcal
{N}(A,1),1\}}\bigr),\nonumber
\end{eqnarray}
obtaining~(\ref{x11}), and the proof is finished.
\end{pf*}

\section*{Acknowledgments}

The authors would like to thank Jiefeng Jiang, two anonymous referees,
the Associate Editor and the Editor for very useful comments.



\printaddresses

\end{document}